\documentclass{gtart}

\def\ifplaintex{\expandafter\ifx\csname documentclass\endcsname\relax}

\def\gtp{{\mathsurround=0pt\it $\cal G\mskip-2mu$eometry \&\ 
$\cal T\!\!$opology $\cal P\!$ublications}}  

\def\recd{{\small Received:\qua\receiveddate\ifx\reviseddate\relax
\else\qquad Revised:\qua\reviseddate\fi\par}} 

\def\lognumber#1{\def\thelognumber{#1}}
\def\volumenumber#1{\def\thevolumenumber{#1}}
\def\volumeyear#1{\def\thevolumeyear{#1}}
\def\papernumber#1{\def\thepapernumber{#1}}
\def\pagenumbers#1#2{\def\startpage{#1}\def\finishpage{#2}}
\def\published#1{\def\publishdate{#1}}

\def\received#1{\def\receiveddate{#1}}
\def\revised#1{\def\reviseddate{#1}}
\def\accepted#1{\def\accepteddate{#1}}

\long\def\asciiabstract#1{\long\def\theasciiabstract{#1}}

\let\\\par\let\thelognumber\relax\let\thevolumenumber\relax
\let\thepapernumber\relax\let\thevolumeyear\relax\let\startpage\relax
\let\finishpage\relax\let\publishdate\relax\let\receiveddate\relax
\let\reviseddate\relax\let\accepteddate\relax\let\theasciititle\relax
\let\theasciiauthors\relax
\let\theasciiabstract\relax

\let\theasciiemail\relax

\ifplaintex
\font\logobig=cmssbx10 scaled 3836
\font\logomed=cmssbx10 scaled 2557
\else
\font\logobig=cmssbx10 scaled 4200
\font\logomed=cmssbx10 scaled 2800
\fi

\long\def\makeagttitle{   
\count0=\startpage
\agt\hfill      
\hbox to 45truept{\vbox to 0pt{\vglue -13truept{\logomed A\kern -.37em{\logobig 
T}\kern -.38em G}\vss}\hss}
\break
{\small Volume \thevolumenumber\ (\thevolumeyear)
\startpage--\finishpage\nl
Published: \publishdate}

\vglue .25truein

{\parskip=0pt\leftskip 0pt plus
1fil\def\\{\par\smallskip}{\Large\bf\thetitle}\par\medskip} \vglue
0.05truein

{\parskip=0pt\leftskip 0pt plus 1fil\def\\{\par}{\sc\theauthors}
\par\medskip}%
 
\vglue 0.03truein 

{\small\leftskip 25truept\rightskip 25truept{\bf Abstract}\stdspace\theabstract

{\bf AMS Classification}\stdspace\theprimaryclass
\ifx\thesecondaryclass\relax\else; \thesecondaryclass\fi\par
{\bf Keywords}\stdspace \thekeywords\par}\vglue 7truept

}   

\ifplaintex
\font\phead=cmsl9 scaled 950
\font\pnum=cmbx10 scaled 913
\font\pfoot=cmsl9 scaled 950
\headline{\vbox to 0pt{\vskip -4.5mm\line{\small\phead\ifnum
\count0=\startpage ISSN 1472-2739 (on-line) 1472-2747 (printed)
\hfill {\pnum\folio}\else\ifodd\count0\def\\{ }%
\ifx\theshorttitle\relax\thetitle\else\theshorttitle\fi\hfill{\pnum\folio}
\else\def\\{ and }{\pnum\folio}\hfill\ifx\theshortauthors\relax\theauthors
\else\theshortauthors\fi\fi\fi}\vss}}
\footline{\vbox to 0pt{\vglue 0mm\line{\small\pfoot\ifnum\count0=\startpage
\copyright\ \gtp\hfill\else
\agt, Volume \thevolumenumber\ (\thevolumeyear)\hfill\fi}\vss}}
\else
\font=cmsl9 scaled 1050
\font\lnum=cmbx10 
\font=cmsl9 scaled 1050
\makeatletter
\def\@oddhead{{\small\lhead\ifnum\count0=\startpage ISSN 1472-2739 
(on-line) 1472-2747 (printed)\hfill {\lnum\number\count0}\else\ifodd\count0
\def\\{ }\ifx\theshorttitle\relax \thetitle \else\theshorttitle\fi\hfill
{\lnum\number\count0}\else\def\\{ and }{\lnum\number\count0}
\hfill\ifx\theshortauthors\relax 
\theauthors\else\theshortauthors\fi\fi\fi}}\def\@evenhead{\@oddhead}
\def\@oddfoot{\small\lfoot\ifnum\count0=\startpage\copyright\ \gtp\hfill\else
\agt, Volume \thevolumenumber\ (\thevolumeyear)\hfill\fi}
\def\@evenfoot{\@oddfoot}
\makeatother
\fi
\let\maketitlepage\makeagttitle
\let\makeshorttitle\maketitlepage
\let\maketitle\maketitlepage

\newwrite\gtoutfile
\long\gdef\makeheadfile{  
{\def\\{, }\def\s{ }
\immediate\openout\gtoutfile head.xxx
\immediate\write\gtoutfile{To: math@arxiv.org}
\immediate\write\gtoutfile{Subject: put OR rep NNNNN:ppppp}
\immediate\write\gtoutfile{--text follows this line--}
\immediate\write\gtoutfile{Proxy-for: \ifx\theasciiauthors\relax
\theauthors\else\theasciiauthors\fi\s<\ifx\theasciiemail\relax\theemail\else\theasciiemail\fi>}
\immediate\write\gtoutfile{\noexpand\\}
\immediate\write\gtoutfile{Authors: \ifx\theasciiauthors\relax
\theauthors\else\theasciiauthors\fi}
{\def\\{ }\immediate\write\gtoutfile{Title: \ifx\theasciititle\relax
\thetitle\else\theasciititle\fi}}
\immediate\write\gtoutfile{Subj-class: GT or SG, GR etc}
\immediate\write\gtoutfile{MSC-class: \theprimaryclass\ifx\thesecondaryclass\relax\else, \thesecondaryclass\fi}
\immediate\write\gtoutfile{Journal-ref: Algebr. Geom. Topol. \thevolumenumber\s
(\thevolumeyear) \startpage-\finishpage}
\immediate\write\gtoutfile{Comments: Published by Algebraic and
Geometric Topology at}
\immediate\write\gtoutfile{\s\s\s  http://www.maths.warwick.ac.uk/agt/AGTVol\thevolumenumber/agt-\thevolumenumber-\thepapernumber.abs.html}
\immediate\write\gtoutfile{\noexpand\\}
\immediate\write\gtoutfile{}
\ifx\theasciiabstract\relax
\immediate\write\gtoutfile{\theabstract}\else
\immediate\write\gtoutfile{\theasciiabstract}\fi
\immediate\write\gtoutfile{}
\immediate\write\gtoutfile{\noexpand\\}
\immediate\write\gtoutfile{}
\immediate\closeout\gtoutfile}}  

\def\maketitlepage{\makeagttitle\makeheadfile}
\let\makeshorttitle\maketitlepage
\let\maketitle\maketitlepage

\def\ifplaintex{\expandafter\ifx\csname documentclass\endcsname\relax}

\def\gtp{{\mathsurround=0pt\it $\cal G\mskip-2mu$eometry \&\ 
$\cal T\!\!$opology $\cal P\!$ublications}}  

\def\recd{{\small Received:\qua\receiveddate\ifx\reviseddate\relax
\else\qquad Revised:\qua\reviseddate\fi\par}} 

\def\lognumber#1{\def\thelognumber{#1}}
\def\volumenumber#1{\def\thevolumenumber{#1}}
\def\volumeyear#1{\def\thevolumeyear{#1}}
\def\papernumber#1{\def\thepapernumber{#1}}
\def\pagenumbers#1#2{\def\startpage{#1}\def\finishpage{#2}}
\def\published#1{\def\publishdate{#1}}

\def\received#1{\def\receiveddate{#1}}
\def\revised#1{\def\reviseddate{#1}}
\def\accepted#1{\def\accepteddate{#1}}

\long\def\asciiabstract#1{\long\def\theasciiabstract{#1}}

\let\\\par\let\thelognumber\relax\let\thevolumenumber\relax
\let\thepapernumber\relax\let\thevolumeyear\relax\let\startpage\relax
\let\finishpage\relax\let\publishdate\relax\let\receiveddate\relax
\let\reviseddate\relax\let\accepteddate\relax\let\theasciititle\relax
\let\theasciiauthors\relax
\let\theasciiabstract\relax

\let\theasciiemail\relax

\ifplaintex
\font\logobig=cmssbx10 scaled 3836
\font\logomed=cmssbx10 scaled 2557
\else
\font\logobig=cmssbx10 scaled 4200
\font\logomed=cmssbx10 scaled 2800
\fi

\long\def\makeagttitle{   
\count0=\startpage
\agt\hfill      
\hbox to 45truept{\vbox to 0pt{\vglue -13truept{\logomed A\kern -.37em{\logobig 
T}\kern -.38em G}\vss}\hss}
\break
{\small Volume \thevolumenumber\ (\thevolumeyear)
\startpage--\finishpage\nl
Published: \publishdate}

\vglue .25truein

{\parskip=0pt\leftskip 0pt plus
1fil\def\\{\par\smallskip}{\Large\bf\thetitle}\par\medskip} \vglue
0.05truein

{\parskip=0pt\leftskip 0pt plus 1fil\def\\{\par}{\sc\theauthors}
\par\medskip}%
 
\vglue 0.03truein 

{\small\leftskip 25truept\rightskip 25truept{\bf Abstract}\stdspace\theabstract

{\bf AMS Classification}\stdspace\theprimaryclass
\ifx\thesecondaryclass\relax\else; \thesecondaryclass\fi\par
{\bf Keywords}\stdspace \thekeywords\par}\vglue 7truept

}   

\ifplaintex
\font\phead=cmsl9 scaled 950
\font\pnum=cmbx10 scaled 913
\font\pfoot=cmsl9 scaled 950
\headline{\vbox to 0pt{\vskip -4.5mm\line{\small\phead\ifnum
\count0=\startpage ISSN 1472-2739 (on-line) 1472-2747 (printed)
\hfill {\pnum\folio}\else\ifodd\count0\def\\{ }%
\ifx\theshorttitle\relax\thetitle\else\theshorttitle\fi\hfill{\pnum\folio}
\else\def\\{ and }{\pnum\folio}\hfill\ifx\theshortauthors\relax\theauthors
\else\theshortauthors\fi\fi\fi}\vss}}
\footline{\vbox to 0pt{\vglue 0mm\line{\small\pfoot\ifnum\count0=\startpage
\copyright\ \gtp\hfill\else
\agt, Volume \thevolumenumber\ (\thevolumeyear)\hfill\fi}\vss}}
\else
\font=cmsl9 scaled 1050
\font\lnum=cmbx10 
\font=cmsl9 scaled 1050
\makeatletter
\def\@oddhead{{\small\lhead\ifnum\count0=\startpage ISSN 1472-2739 
(on-line) 1472-2747 (printed)\hfill {\lnum\number\count0}\else\ifodd\count0
\def\\{ }\ifx\theshorttitle\relax \thetitle \else\theshorttitle\fi\hfill
{\lnum\number\count0}\else\def\\{ and }{\lnum\number\count0}
\hfill\ifx\theshortauthors\relax 
\theauthors\else\theshortauthors\fi\fi\fi}}\def\@evenhead{\@oddhead}
\def\@oddfoot{\small\lfoot\ifnum\count0=\startpage\copyright\ \gtp\hfill\else
\agt, Volume \thevolumenumber\ (\thevolumeyear)\hfill\fi}
\def\@evenfoot{\@oddfoot}
\makeatother
\fi
\let\maketitlepage\makeagttitle
\let\makeshorttitle\maketitlepage
\let\maketitle\maketitlepage

\newwrite\gtoutfile
\long\gdef\makeheadfile{  
{\def\\{, }\def\s{ }
\immediate\openout\gtoutfile head.xxx
\immediate\write\gtoutfile{To: math@arxiv.org}
\immediate\write\gtoutfile{Subject: put OR rep NNNNN:ppppp}
\immediate\write\gtoutfile{--text follows this line--}
\immediate\write\gtoutfile{Proxy-for: \ifx\theasciiauthors\relax
\theauthors\else\theasciiauthors\fi\s<\ifx\theasciiemail\relax\theemail\else\theasciiemail\fi>}
\immediate\write\gtoutfile{\noexpand\\}
\immediate\write\gtoutfile{Authors: \ifx\theasciiauthors\relax
\theauthors\else\theasciiauthors\fi}
{\def\\{ }\immediate\write\gtoutfile{Title: \ifx\theasciititle\relax
\thetitle\else\theasciititle\fi}}
\immediate\write\gtoutfile{Subj-class: GT or SG, GR etc}
\immediate\write\gtoutfile{MSC-class: \theprimaryclass\ifx\thesecondaryclass\relax\else, \thesecondaryclass\fi}
\immediate\write\gtoutfile{Journal-ref: Algebr. Geom. Topol. \thevolumenumber\s
(\thevolumeyear) \startpage-\finishpage}
\immediate\write\gtoutfile{Comments: Published by Algebraic and
Geometric Topology at}
\immediate\write\gtoutfile{\s\s\s  http://www.maths.warwick.ac.uk/agt/AGTVol\thevolumenumber/agt-\thevolumenumber-\thepapernumber.abs.html}
\immediate\write\gtoutfile{\noexpand\\}
\immediate\write\gtoutfile{}
\ifx\theasciiabstract\relax
\immediate\write\gtoutfile{\theabstract}\else
\immediate\write\gtoutfile{\theasciiabstract}\fi
\immediate\write\gtoutfile{}
\immediate\write\gtoutfile{\noexpand\\}
\immediate\write\gtoutfile{}
\immediate\closeout\gtoutfile}}  

\def\maketitlepage{\makeagttitle\makeheadfile}
\let\makeshorttitle\maketitlepage
\let\maketitle\maketitlepage

\def\ifplaintex{\expandafter\ifx\csname documentclass\endcsname\relax}

\def\gtp{{\mathsurround=0pt\it $\cal G\mskip-2mu$eometry \&\ 
$\cal T\!\!$opology $\cal P\!$ublications}}  

\def\recd{{\small Received:\qua\receiveddate\ifx\reviseddate\relax
\else\qquad Revised:\qua\reviseddate\fi\par}} 

\def\lognumber#1{\def\thelognumber{#1}}
\def\volumenumber#1{\def\thevolumenumber{#1}}
\def\volumeyear#1{\def\thevolumeyear{#1}}
\def\papernumber#1{\def\thepapernumber{#1}}
\def\pagenumbers#1#2{\def\startpage{#1}\def\finishpage{#2}}
\def\published#1{\def\publishdate{#1}}

\def\received#1{\def\receiveddate{#1}}
\def\revised#1{\def\reviseddate{#1}}
\def\accepted#1{\def\accepteddate{#1}}

\long\def\asciiabstract#1{\long\def\theasciiabstract{#1}}

\let\\\par\let\thelognumber\relax\let\thevolumenumber\relax
\let\thepapernumber\relax\let\thevolumeyear\relax\let\startpage\relax
\let\finishpage\relax\let\publishdate\relax\let\receiveddate\relax
\let\reviseddate\relax\let\accepteddate\relax\let\theasciititle\relax
\let\theasciiauthors\relax
\let\theasciiabstract\relax

\let\theasciiemail\relax

\ifplaintex
\font\logobig=cmssbx10 scaled 3836
\font\logomed=cmssbx10 scaled 2557
\else
\font\logobig=cmssbx10 scaled 4200
\font\logomed=cmssbx10 scaled 2800
\fi

\long\def\makeagttitle{   
\count0=\startpage
\agt\hfill      
\hbox to 45truept{\vbox to 0pt{\vglue -13truept{\logomed A\kern -.37em{\logobig 
T}\kern -.38em G}\vss}\hss}
\break
{\small Volume \thevolumenumber\ (\thevolumeyear)
\startpage--\finishpage\nl
Published: \publishdate}

\vglue .25truein

{\parskip=0pt\leftskip 0pt plus
1fil\def\\{\par\smallskip}{\Large\bf\thetitle}\par\medskip} \vglue
0.05truein

{\parskip=0pt\leftskip 0pt plus 1fil\def\\{\par}{\sc\theauthors}
\par\medskip}%
 
\vglue 0.03truein 

{\small\leftskip 25truept\rightskip 25truept{\bf Abstract}\stdspace\theabstract

{\bf AMS Classification}\stdspace\theprimaryclass
\ifx\thesecondaryclass\relax\else; \thesecondaryclass\fi\par
{\bf Keywords}\stdspace \thekeywords\par}\vglue 7truept

}   

\ifplaintex
\font\phead=cmsl9 scaled 950
\font\pnum=cmbx10 scaled 913
\font\pfoot=cmsl9 scaled 950
\headline{\vbox to 0pt{\vskip -4.5mm\line{\small\phead\ifnum
\count0=\startpage ISSN 1472-2739 (on-line) 1472-2747 (printed)
\hfill {\pnum\folio}\else\ifodd\count0\def\\{ }%
\ifx\theshorttitle\relax\thetitle\else\theshorttitle\fi\hfill{\pnum\folio}
\else\def\\{ and }{\pnum\folio}\hfill\ifx\theshortauthors\relax\theauthors
\else\theshortauthors\fi\fi\fi}\vss}}
\footline{\vbox to 0pt{\vglue 0mm\line{\small\pfoot\ifnum\count0=\startpage
\copyright\ \gtp\hfill\else
\agt, Volume \thevolumenumber\ (\thevolumeyear)\hfill\fi}\vss}}
\else
\font=cmsl9 scaled 1050
\font\lnum=cmbx10 
\font=cmsl9 scaled 1050
\makeatletter
\def\@oddhead{{\small\lhead\ifnum\count0=\startpage ISSN 1472-2739 
(on-line) 1472-2747 (printed)\hfill {\lnum\number\count0}\else\ifodd\count0
\def\\{ }\ifx\theshorttitle\relax \thetitle \else\theshorttitle\fi\hfill
{\lnum\number\count0}\else\def\\{ and }{\lnum\number\count0}
\hfill\ifx\theshortauthors\relax 
\theauthors\else\theshortauthors\fi\fi\fi}}\def\@evenhead{\@oddhead}
\def\@oddfoot{\small\lfoot\ifnum\count0=\startpage\copyright\ \gtp\hfill\else
\agt, Volume \thevolumenumber\ (\thevolumeyear)\hfill\fi}
\def\@evenfoot{\@oddfoot}
\makeatother
\fi
\let\maketitlepage\makeagttitle
\let\makeshorttitle\maketitlepage
\let\maketitle\maketitlepage

\newwrite\gtoutfile
\long\gdef\makeheadfile{  
{\def\\{, }\def\s{ }
\immediate\openout\gtoutfile head.xxx
\immediate\write\gtoutfile{To: math@arxiv.org}
\immediate\write\gtoutfile{Subject: put OR rep NNNNN:ppppp}
\immediate\write\gtoutfile{--text follows this line--}
\immediate\write\gtoutfile{Proxy-for: \ifx\theasciiauthors\relax
\theauthors\else\theasciiauthors\fi\s<\ifx\theasciiemail\relax\theemail\else\theasciiemail\fi>}
\immediate\write\gtoutfile{\noexpand\\}
\immediate\write\gtoutfile{Authors: \ifx\theasciiauthors\relax
\theauthors\else\theasciiauthors\fi}
{\def\\{ }\immediate\write\gtoutfile{Title: \ifx\theasciititle\relax
\thetitle\else\theasciititle\fi}}
\immediate\write\gtoutfile{Subj-class: GT or SG, GR etc}
\immediate\write\gtoutfile{MSC-class: \theprimaryclass\ifx\thesecondaryclass\relax\else, \thesecondaryclass\fi}
\immediate\write\gtoutfile{Journal-ref: Algebr. Geom. Topol. \thevolumenumber\s
(\thevolumeyear) \startpage-\finishpage}
\immediate\write\gtoutfile{Comments: Published by Algebraic and
Geometric Topology at}
\immediate\write\gtoutfile{\s\s\s  http://www.maths.warwick.ac.uk/agt/AGTVol\thevolumenumber/agt-\thevolumenumber-\thepapernumber.abs.html}
\immediate\write\gtoutfile{\noexpand\\}
\immediate\write\gtoutfile{}
\ifx\theasciiabstract\relax
\immediate\write\gtoutfile{\theabstract}\else
\immediate\write\gtoutfile{\theasciiabstract}\fi
\immediate\write\gtoutfile{}
\immediate\write\gtoutfile{\noexpand\\}
\immediate\write\gtoutfile{}
\immediate\closeout\gtoutfile}}  

\def\maketitlepage{\makeagttitle\makeheadfile}
\let\makeshorttitle\maketitlepage
\let\maketitle\maketitlepage

\def\ifplaintex{\expandafter\ifx\csname documentclass\endcsname\relax}

\def\gtp{{\mathsurround=0pt\it $\cal G\mskip-2mu$eometry \&\ 
$\cal T\!\!$opology $\cal P\!$ublications}}  

\def\recd{{\small Received:\qua\receiveddate\ifx\reviseddate\relax
\else\qquad Revised:\qua\reviseddate\fi\par}} 

\def\lognumber#1{\def\thelognumber{#1}}
\def\volumenumber#1{\def\thevolumenumber{#1}}
\def\volumeyear#1{\def\thevolumeyear{#1}}
\def\papernumber#1{\def\thepapernumber{#1}}
\def\pagenumbers#1#2{\def\startpage{#1}\def\finishpage{#2}}
\def\published#1{\def\publishdate{#1}}

\def\received#1{\def\receiveddate{#1}}
\def\revised#1{\def\reviseddate{#1}}
\def\accepted#1{\def\accepteddate{#1}}

\long\def\asciiabstract#1{\long\def\theasciiabstract{#1}}

\let\\\par\let\thelognumber\relax\let\thevolumenumber\relax
\let\thepapernumber\relax\let\thevolumeyear\relax\let\startpage\relax
\let\finishpage\relax\let\publishdate\relax\let\receiveddate\relax
\let\reviseddate\relax\let\accepteddate\relax\let\theasciititle\relax
\let\theasciiauthors\relax
\let\theasciiabstract\relax

\let\theasciiemail\relax

\ifplaintex
\font\logobig=cmssbx10 scaled 3836
\font\logomed=cmssbx10 scaled 2557
\else
\font\logobig=cmssbx10 scaled 4200
\font\logomed=cmssbx10 scaled 2800
\fi

\long\def\makeagttitle{   
\count0=\startpage
\agt\hfill      
\hbox to 45truept{\vbox to 0pt{\vglue -13truept{\logomed A\kern -.37em{\logobig 
T}\kern -.38em G}\vss}\hss}
\break
{\small Volume \thevolumenumber\ (\thevolumeyear)
\startpage--\finishpage\nl
Published: \publishdate}

\vglue .25truein

{\parskip=0pt\leftskip 0pt plus
1fil\def\\{\par\smallskip}{\Large\bf\thetitle}\par\medskip} \vglue
0.05truein

{\parskip=0pt\leftskip 0pt plus 1fil\def\\{\par}{\sc\theauthors}
\par\medskip}%
 
\vglue 0.03truein 

{\small\leftskip 25truept\rightskip 25truept{\bf Abstract}\stdspace\theabstract

{\bf AMS Classification}\stdspace\theprimaryclass
\ifx\thesecondaryclass\relax\else; \thesecondaryclass\fi\par
{\bf Keywords}\stdspace \thekeywords\par}\vglue 7truept

}   

\ifplaintex
\font\phead=cmsl9 scaled 950
\font\pnum=cmbx10 scaled 913
\font\pfoot=cmsl9 scaled 950
\headline{\vbox to 0pt{\vskip -4.5mm\line{\small\phead\ifnum
\count0=\startpage ISSN 1472-2739 (on-line) 1472-2747 (printed)
\hfill {\pnum\folio}\else\ifodd\count0\def\\{ }%
\ifx\theshorttitle\relax\thetitle\else\theshorttitle\fi\hfill{\pnum\folio}
\else\def\\{ and }{\pnum\folio}\hfill\ifx\theshortauthors\relax\theauthors
\else\theshortauthors\fi\fi\fi}\vss}}
\footline{\vbox to 0pt{\vglue 0mm\line{\small\pfoot\ifnum\count0=\startpage
\copyright\ \gtp\hfill\else
\agt, Volume \thevolumenumber\ (\thevolumeyear)\hfill\fi}\vss}}
\else
\font=cmsl9 scaled 1050
\font\lnum=cmbx10 
\font=cmsl9 scaled 1050
\makeatletter
\def\@oddhead{{\small\lhead\ifnum\count0=\startpage ISSN 1472-2739 
(on-line) 1472-2747 (printed)\hfill {\lnum\number\count0}\else\ifodd\count0
\def\\{ }\ifx\theshorttitle\relax \thetitle \else\theshorttitle\fi\hfill
{\lnum\number\count0}\else\def\\{ and }{\lnum\number\count0}
\hfill\ifx\theshortauthors\relax 
\theauthors\else\theshortauthors\fi\fi\fi}}\def\@evenhead{\@oddhead}
\def\@oddfoot{\small\lfoot\ifnum\count0=\startpage\copyright\ \gtp\hfill\else
\agt, Volume \thevolumenumber\ (\thevolumeyear)\hfill\fi}
\def\@evenfoot{\@oddfoot}
\makeatother
\fi
\let\maketitlepage\makeagttitle
\let\makeshorttitle\maketitlepage
\let\maketitle\maketitlepage

\newwrite\gtoutfile
\long\gdef\makeheadfile{  
{\def\\{, }\def\s{ }
\immediate\openout\gtoutfile head.xxx
\immediate\write\gtoutfile{To: math@arxiv.org}
\immediate\write\gtoutfile{Subject: put OR rep NNNNN:ppppp}
\immediate\write\gtoutfile{--text follows this line--}
\immediate\write\gtoutfile{Proxy-for: \ifx\theasciiauthors\relax
\theauthors\else\theasciiauthors\fi\s<\ifx\theasciiemail\relax\theemail\else\theasciiemail\fi>}
\immediate\write\gtoutfile{\noexpand\\}
\immediate\write\gtoutfile{Authors: \ifx\theasciiauthors\relax
\theauthors\else\theasciiauthors\fi}
{\def\\{ }\immediate\write\gtoutfile{Title: \ifx\theasciititle\relax
\thetitle\else\theasciititle\fi}}
\immediate\write\gtoutfile{Subj-class: GT or SG, GR etc}
\immediate\write\gtoutfile{MSC-class: \theprimaryclass\ifx\thesecondaryclass\relax\else, \thesecondaryclass\fi}
\immediate\write\gtoutfile{Journal-ref: Algebr. Geom. Topol. \thevolumenumber\s
(\thevolumeyear) \startpage-\finishpage}
\immediate\write\gtoutfile{Comments: Published by Algebraic and
Geometric Topology at}
\immediate\write\gtoutfile{\s\s\s  http://www.maths.warwick.ac.uk/agt/AGTVol\thevolumenumber/agt-\thevolumenumber-\thepapernumber.abs.html}
\immediate\write\gtoutfile{\noexpand\\}
\immediate\write\gtoutfile{}
\ifx\theasciiabstract\relax
\immediate\write\gtoutfile{\theabstract}\else
\immediate\write\gtoutfile{\theasciiabstract}\fi
\immediate\write\gtoutfile{}
\immediate\write\gtoutfile{\noexpand\\}
\immediate\write\gtoutfile{}
\immediate\closeout\gtoutfile}}  

\def\maketitlepage{\makeagttitle\makeheadfile}
\let\makeshorttitle\maketitlepage
\let\maketitle\maketitlepage

\lognumber{26}
\volumenumber{2}
\volumeyear{2002}
\papernumber{26}
\published{19 July 2002}
\pagenumbers{537}{562}
\received{3 October 2001}
\revised{29 April 2002}
\accepted{26 June 2002}

\usepackage {amsmath,amssymb}

\newtheorem{lemma}{Lemma}[section]
\newtheorem{theorem}[lemma]{Theorem}
\newtheorem{proposition}[lemma]{Proposition}
\newtheorem{corollary}[lemma]{Corollary}

\theoremstyle{definition}

\newtheorem*{theorem*}{Theorem}
\newtheorem*{remark}{Remark}

\newcommand{\romanIndex}
 {\renewcommand{\labelenumi}{\rm(\roman{enumi})}}

\begin{document}
\title{Farrell cohomology of low genus pure 
mapping\\class groups with punctures}
\shorttitle{Farrell cohomology of mapping class groups}
\author {Qin Lu}
\address{Department of Mathematics, Lafayette College\\Easton, PA 18042, USA}
\email{luq@Lafayette.edu}

\begin{abstract} 
In this paper, we calculate the $p$-torsion of the Farrell cohomology
for low genus pure mapping class groups with punctures, where $p$ is
an odd prime. Here, `low genus' means $g=1,2,3$; and `pure mapping
class groups with punctures' means the mapping class groups with any
number of punctures, where the punctures are not allowed to be
permuted. These calculations use our previous results about the
periodicity of pure mapping class groups with punctures, as well as
other cohomological tools. The low genus cases are interesting because
we know that the high genus cases can be reduced to the low genus
ones. Also, the cohomological properties of the mapping class groups
without punctures are closely related to our cases.
\end{abstract}

\asciiabstract{ 
In this paper, we calculate the p-torsion of the Farrell cohomology
for low genus pure mapping class groups with punctures, where p is an
odd prime. Here, `low genus' means g=1,2,3; and `pure mapping class
groups with punctures' means the mapping class groups with any number
of punctures, where the punctures are not allowed to be
permuted. These calculations use our previous results about the
periodicity of pure mapping class groups with punctures, as well as
other cohomological tools. The low genus cases are interesting because
we know that the high genus cases can be reduced to the low genus
ones. Also, the cohomological properties of the mapping class groups
without punctures are closely related to our cases.}



\maketitle

\section*{Introduction}

The pure mapping class group with punctures, $\Gamma^i_g,\ \ $ is defined as 
$$\pi _0 (\hbox{\it Diffeo}^+(S_g,P_1,P_2,...P_i)),$$ 
where $\hbox{\it Diffeo}^+(S_g,P_1,P_2,...P_i)$ is the group of orientation preserving diffeomorphisms of $S_g$ (closed orientable two manifold with genus $g$) which fix
the points $P_j$ individually. For $i\ge 1$, we refer to $\Gamma_g^i$ as the $pure$ mapping class group with $punctures$. We write $\Gamma_g=\Gamma_g^0$, which we refer to as the $unpunctured$ mapping class group. We also write $\tilde \Gamma^i_g$ as the mapping class group with $punctures$, where the punctures are allowed to be permuted.

Recall that a group $\Gamma$ of finite virtual cohomological dimension
is said to be periodic (in cohomology) if for some $d\neq 0$ there is
an element $u\in \hat H^d(\Gamma, \mathbb Z)$ which is invertible in
the ring $\hat H^*(\Gamma,\mathbb Z)$. Cup product with $u$ then gives
a periodicity isomorphism $$\hat H^i(\Gamma, M)\cong \hat
H^{i+d}(\Gamma, M)$$ for any $\Gamma-$module $M$ and any $i\in \mathbb
Z.$ Similarly, we say that $\Gamma$ is $p$-periodic (where $p$ is a
prime) if the $p$-primary component $\hat H^*(\Gamma, \mathbb
Z)_{(p)},$ which is itself a ring, contains an invertible element of
non-zero degree $d$. We then have $$\hat H^i(\Gamma, M)_{(p)}\cong
\hat H^{i+d}(\Gamma, M)_{(p)}.$$ We refer to $d$ as the period(the
$p$-period) of the group $\Gamma$.

It is known that $\Gamma_g$ is never $2$-periodic for $g>0$. For an odd prime $p$, $\Gamma_g$ is
 $p$-periodic if and only if $g$ and $p$ satisfy certain relations. Moreover, 
the $p$-period depends on the genus $g$. However, we proved that the pure mapping class group with punctures is periodic and the period is $2$. [4]. Because of this property, it is only necessary to calculate a small range of cohomology groups, namely the even dimension and the odd dimension of cohomology groups, in order to determine the entire additive cohomology structure; this makes the calculation of $p$-torsion of the Farrell cohomology of the pure mapping class group with punctures possible. However, calculating cohomology is very hard in general. We will only calculate the low genus cases with $p$ odd. In fact, high genus cases can be reduced to low genus cases so our results can be further generalized. The case $p=2$ is very complicated, and we will 
not calculate it. 

In this paper, we calculate Farrell cohomology. This agrees with the standard cohomology above the finite virtual cohomological dimension(vcd). It is well known that any mapping class group has finite vcd and the vcd has been calculated explicitly. 
 
We will use the following theorem in K.S.Brown's book [1].

If $\Gamma$ is a $p$-periodic group, then 
$$\hat H^*(\Gamma)_{(p)}\buildrel \cong \over \longrightarrow \Pi_{\mathbb Z/p\in S}\hat H^*(N(\mathbb Z/p))_{(p)},$$
where $\hat H^*(\Gamma)_{(p)}$ stands for the $p$-torsion of the Farrell cohomology of $\Gamma$, $S$ is a set of representatives for the conjugacy classes of subgroups of $\Gamma$ of order
$p$, and $N(\mathbb Z/p)$ is the normalizer of $\mathbb Z/p$ in $\Gamma$.  

By [4], we know that $\Gamma_g^i$($i\ge 1$) is periodic, thus $p$-periodic for any prime $p$. Hence, we can apply the above theorem to our calculation.
Of course, one must be careful with the details.

The paper is divided into two sections. In the first, we analyze the
$p$-torsion in $\Gamma_g^i\ \ (i\ge 1,g\ge 1),$ where $p$ is any
prime. The basic tools are the Riemann Hurwitz Equation [3], Nielsen's
Realization Theorem and some results in [4] related to the pure
mapping class groups with punctures. In section 2, we calculate the
$p$-torsion of the Farrell cohomology of $\Gamma_g^i,~i\ge1,~g=1,2,3.$
For this we need to analyze the properties of the normalizer of the
subgroup of order $p$ in $\Gamma_g^i$. A result of MacLachlan and
Harvey [5] states that for $\mathbb Z/p<\Gamma_g^i$, the quotient
$N(\mathbb Z/p)/(\mathbb Z/p)$ maps injectively into the mapping class
group $\tilde \Gamma_h^t$, where $h$ is the genus of the orbit space
$S_g/(\mathbb Z/p)$, and $t$ is the number of fixed points. Note that
$h\le g$, so the high genus cases can be reduced to the low genus
cases. Using the properties of $\tilde \Gamma_h^t$ which we develop
later, we find $\hat H^*(N(\mathbb Z/p)/(\mathbb Z/p), F_p)$. Then, by
the short exact sequence $1\to \mathbb Z/p\to N(\mathbb Z/p)\to
N(\mathbb Z/p)/(\mathbb Z/p)\to 1$, we can calculate $\hat
H^*(N(\mathbb Z/p), \mathbb Z)_{(p)}$. The basic tools here are
Cohen's and Xia's results for mapping class groups [2, 7], cohomology
of symmetric groups, and the Serre spectral sequence. In order to
finish our calculation, we also need to count the number of conjugacy
classes of subgroups of $\Gamma_g^i$ of order $p$. This is related to
the fixed point data of the pure mapping class group with punctures
[4].

Fixed point data have been well-defined in [6] for the unpunctured mapping class group. In [4], we generalized the fixed point data to the case of pure mapping class group with punctures. Recall that for an element of order $p$, $\alpha \in \Gamma_g^i,$ we can lift $\alpha$ to $f$, an orientation-preserving 
diffeomorphism of the closed orientable surface $S_g$ of prime period $p$. Note that by the definition of $\Gamma_g^i$, $f$ has already fixed $i$ points. Assume that $f$ acts on $S_g$ with $t$ fixed points total. The fixed point data of $f$ are a set $\delta (f)=(\beta_1,..., \beta_i|\beta_{i+1},...\beta_{t})$, where $t$ is the number of fixed points of $f$; $\beta_1,..., \beta_i$ are ordered, corresponding to the $i$ fixed points associated to $\Gamma_g^i$;  $\beta_{i+1},...\beta_{t}$ are unordered, corresponding to the rest $t-i$ fixed points which the $f$-action on $S_g$ has. Each $\beta_j$ is the integer (mod $p$) such that $f^{\beta_j}$ acts as multiplication by $e^{2\pi i/p}$ in the local invariant complex structure at the jth fixed point. In [4], we proved that the fixed point data are well defined for $\alpha \in \Gamma_g^i,$ which is induced by the fixed point data of $f$. Moreover, for any subgroup of order $p$, we can pick a generator $\alpha$, such that $\delta (\alpha)=(1,\beta_2,..., \beta_i|\beta_{i+1},...\beta_{t}),$ namely $\beta_1=1$. From now on, we may assume $\beta_1=1$ for our fixed point data. 

By Theorem 2.5 and Proposition 2.6 in [4], we can  count the  
conjugacy classes of subgroups of $\Gamma_g^i$ of order $p$ by using the fixed point data.

\section{\bf{The $p$-torsion in $\Gamma_g^i(i\ge 1,g>0).$}}

In this section, we investigate the $p$-torsion in $\Gamma_1^i$, $\Gamma_2^i$, 
$\Gamma_3^i$ for $i\ge 1$. The basic tools are the Riemann Hurwitz Equation and Nielsen's Realization Theorem.

\begin{lemma}
\label{a} 
\begin{enumerate}\romanIndex
\item If $\Gamma_g^i(i> 2)$ has $p$-torsion, 

then $p\le 2g/(i-2)+1$, where $p$ is any prime and $g>0.$
\item If $\Gamma_g^i(i=1, i=2)$ has $p$-torsion, 

then $p\le 2g+1$, where $p$ is any prime and $g>0.$
\end{enumerate}
\end{lemma}
\begin{proof}
(i)\qua If $\Gamma_g^i$ has $p$-torsion, we know that there is $\mathbb Z/p<\Gamma_g^i$.
By Nielsen's Realization theorem [5], we can lift the $\mathbb Z/p< \Gamma_g^i$ 

into $\mathbb Z/p < \hbox{\it Diffeo}^+(S_g,P_1,P_2,...P_i)$. 
Then we can view $\mathbb Z/p$ acting on $S_g$ 
with at least $i$ points fixed. By the property of Riemann Surfaces, the Riemann Hurwitz equation 
$2g-2=p(2h-2)+t(p-1)$ should have positive solutions $(h,t)$, 
where $h$ corresponds to the genus of the quotient space by the $\mathbb Z/p$ action on $S_g$, 
and $t$ is the number of fixed points of this action. 

Since $h\ge 0$ in the Riemann Hurwitz equation, we know that 
$$2g-2\ge -2p+tp-t,$$ 
$$2g+t-2\ge (t-2)p,$$
Since $t$ is the number of fixed points, $t\ge i>2.$
Hence, $$2g/(i-2)+1\ge 2g/(t-2)+1\ge p$$
i.e. $$p\le 2g/(i-2)+1.$$

(ii)\qua By the same argument as in (i), we know that $2g+t-2\ge (t-2)p$.
Since $i=1$ or $i=2$, we have $t\ge 1$. By Theorem 2.7 in [4], we know that if $\Gamma_g^i$ contains a subgroup of order $p$, then the number of fixed points $t$ can not be 1. Hence we only need to consider $t\ge 2$.  
If $t=2$, then by the Riemann Hurwitz equation, $2g-2=p(2h-2)+2(p-1)$ implies $g=ph$. So, $p\le g\le 2g+1$ for $g>0$.
If $t>2$, then $p\le 2g/(t-2)+1$ implies $p\le 2g+1$. 
\end{proof}

\begin{remark}
$H^*(\Gamma_0^n, \mathbb Z)$ is completely calculated by Cohen in [2]. So, in this paper, we
will not consider the case $g=0$. 
The following corollaries determine the $p$-torsion in 
$\Gamma_1^i$, $\Gamma_2^i$, $\Gamma_3^i$ for $i\ge 1$.
\end{remark}

We need to use Theorem 2.7 in [4] for the following corollaries. Hence, we cite it here as a reference.

\begin{theorem}[Theorem 2.7 in \cite{4}]
The Riemann Hurwitz equation $2g-2=p(2h-2)+t(p-1)$ has a non-negative integer solution $(h,t)$, with $t\neq 1$ and $t\ge i$ iff $\Gamma_g^i$ contains a subgroup of order $p$, the subgroup of order $p$ acts on $S_g$ with $t$ fixed points. 
\end{theorem}

\begin{corollary}
\label{b} 
\begin{enumerate}\romanIndex 
\item If $\Gamma_2^i$ has $p$-torsion, then $p=2,3,5$.
\item $\Gamma_2^1$ has 2,3,5 torsion.
\item $\Gamma_2^2$ has 2,3,5 torsion.
\item $\Gamma_2^3$ has 2,3,5 torsion.
\item $\Gamma_2^4$ has 2,3 torsion.
\item $\Gamma_2^5$ has 2 torsion.
\item $\Gamma_2^6$ has 2 torsion.
\item $\Gamma_2^i$ has no $p$-torsion for $i\ge 7$.
\end{enumerate}
\end{corollary}
\begin{proof}

(i)\qua By Lemma \ref{a}, If $i>2$, then $p\le 2g/(i-2)+1=4/(i-2)+1\le 5$; If $i=1,2$, then $p\le 2g+1=5$.

(ii)-(viii)\qua In fact, we are not only interested in the $p$-torsion that $\Gamma_2^i$ contains. We are also interested in 
the values of $t$ and $h$ related to the $\mathbb Z/p$ action, namely, the number of fixed points of the $\mathbb Z/p$ action, and the genus of the quotient space of the $\mathbb Z/p$ action. In (i), we have proved that $\Gamma_2^i$ may contain $2,3,5$ torsion. Now we need to investigate what torsion it indeed contains. We will use the theorem mentioned above.(Theorem 2.7 in [4]) It gives the necessary and sufficient conditions for $\Gamma_g^i$ containing $p$-torsion. We will do this case by case.

Case (1): $p=2$

Plug $g=2$ and $p=2$ into the Riemann Hurwitz equation. Then $2\times 2-2=2(2h-2)+t(2-1)$ implies
$2=4h-4+t,$ i.e., $6=4h+t$. The non-negative integer solutions are $(h,t)=(1,2)$ or $(h,t)=(0,6)$. By Theorem 2.7 in [4], $\Gamma_2^i$ has $2$-torsion for $i\le 6$ and $\Gamma_2^i$ can not have $2$-torsion for $i\ge 7$.

For $\Gamma_2^1$, $\Gamma_2^2$, we have $(h,t)=(1,2)$ or $(0,6)$, so the $\mathbb Z/2$ action on $S_2$ must have $2$ fixed points or $6$ fixed points.
For $\Gamma_2^3$, $\Gamma_2^4$, $\Gamma_2^5$, $\Gamma_2^6$, $(h,t)=(0,6)$, so the $\mathbb Z/2$ action on $S_2$ has $6$ fixed points.

Case (2): $p=3$

As in Case (1), the Riemann Hurwitz equation has non-negative integer solutions: $(h,t)=(1,1)$ or $(h,t)=(0,4)$. By Theorem 2.7 in [4], the only solution will be $(h,t)=(0,4)$.
$\Gamma_2^1$, $\Gamma_2^2$, $\Gamma_2^3$, $\Gamma_2^4$ have $3$-torsion and $\Gamma_2^i$ can not have $3$-torsion for $i\ge 5$. The $\mathbb Z/3$ action on $S_2$ must have $4$ fixed points.

Case (3): $p=5$

Here, the Riemann Hurwitz equation has non-negative integer solutions: 

$(h,t)=(0,3)$. $\Gamma_2^1$, $\Gamma_2^2$, $\Gamma_2^3$ have $5$-torsion and $\Gamma_2^i$ can not have $5$-torsion for $i\ge 4$. 

The $\mathbb Z/5$ action on $S_2$ has $3$ fixed points.
\end{proof}

\begin{corollary}
\label{c} 
\begin{enumerate}\romanIndex 
\item If $\Gamma_3^i$ has $p$-torsion, then $p=2,3,5,7$.
\item $\Gamma_3^1$ has 2,3,7 torsion.
\item $\Gamma_3^2$ has 2,3,7 torsion.
\item $\Gamma_3^3$ has 2,3,7 torsion.
\item $\Gamma_3^4$ has 2,3 torsion.
\item $\Gamma_3^5$ has 2,3 torsion.
\item $\Gamma_3^6$ has 2 torsion.
\item $\Gamma_3^7$ has 2 torsion.
\item $\Gamma_3^8$ has 2 torsion.
\item $\Gamma_3^i$ has no $p$-torsion for $i\ge 9$.
\end{enumerate}
\end{corollary}
\begin{proof}

(i)\qua By Lemma \ref{a}, if $i>2$, then $p\le 2g/(i-2)+1=6/(i-2)+1\le 7$; if $i=1,2$, then $p\le 2g+1=7$.

(ii)-(x)\qua As before, we are not only interested in the $\Gamma_3^i$'s $p$-torsion, but are also interested in 
the value of $t$ and $h$ related to the $\mathbb Z/p$ action. 

Case (1): $p=2$

The Riemann Hurwitz equation has non-negative integer solutions: $(h,t)=(1,4)$ or $(h,t)=(0,8)$. $\Gamma_3^1$, $\Gamma_3^2$, $\Gamma_3^3$, $\Gamma_3^4$, $\Gamma_3^5$, $\Gamma_3^6$, $\Gamma_3^7$, $\Gamma_3^8$ have $2$-torsion and $\Gamma_3^i$ can not have $2$-torsion for $i\ge 9$. 

For $\Gamma_3^1$, $\Gamma_3^2$, $\Gamma_3^3$, $\Gamma_3^4$, $(h,t)=(1,4)$ or $(0,8)$. The $\mathbb Z/2$ action on $S_3$ can have $4$ fixed points or $8$ fixed points.
For $\Gamma_3^5$, $\Gamma_3^6$, $\Gamma_3^7$, $\Gamma_3^8$, $(h,t)=(0,8)$. The $\mathbb Z/2$ action on $S_3$ has $8$ fixed points.

Case (2): $p=3$

The Riemann Hurwitz equation has non-negative integer solutions: $(h,t)=(1,2)$ or $(h,t)=(0,5)$. $\Gamma_3^1$, $\Gamma_3^2$, $\Gamma_3^3$, $\Gamma_3^4$, $\Gamma_3^5$ have $3$-torsion and $\Gamma_3^i$ can not have $3$-torsion for $i\ge 6$. 

For $\Gamma_3^1$, $\Gamma_3^2$, $(h,t)=(1,2)$ or $(0,5)$. The $\mathbb Z/3$ action on $S_3$ can  have $2$ fixed points or $5$ fixed points.
For $\Gamma_3^3$, $\Gamma_3^4$, $\Gamma_3^5$, $(h,t)=(0,5)$. The $\mathbb Z/3$ action on $S_3$ must have $5$ fixed points.

Case (3): $p=5$

The Riemann Hurwitz equation has non-negative integer solution: $(h,t)=(1,1)$. Moreover, as  $t\neq 1$, $\Gamma_3^i$ can not have $5$-torsion.

Case (4): $p=7$

The Riemann Hurwitz equation has non-negative integer solution: $(h,t)=(0,3)$. $\Gamma_3^1$, $\Gamma_3^2$, $\Gamma_3^3$ have $7$-torsion and $\Gamma_3^i$ can not have $7$-torsion for $i\ge 4$. The $\mathbb Z/7$ action on $S_3$ must have $3$ fixed points.
\end{proof}

\begin{corollary}
\label{d} 
\begin{enumerate}\romanIndex 
\item If $\Gamma_1^i$ has $p$-torsion, then $p=2,3$.
\item $\Gamma_1^1$ has 2,3 torsion.
\item $\Gamma_1^2$ has 2,3 torsion.
\item $\Gamma_1^3$ has 2,3 torsion.
\item $\Gamma_1^4$ has 2 torsion.
\item $\Gamma_1^5$ has no $p$-torsion for $i\ge 5$.
\end{enumerate}
\end{corollary}
\begin{proof}
This follows by the same arguments as Corollary \ref{b}.
\end{proof}

\begin{remark}
We summarize the above results. Note that in all our cases $(h,t)=(0,t)$ or $(h,t)=(1,t),$ where $t$ differs case by case.\bigskip

{\small\leftskip 25pt
                  \cl{\bf Table of the solutions of $(h,t)$ for Riemann Hurwitz equation}

\ \ \ \ \ \ \ \ \ \ \ \ \ \ \ \ \ \ \ \ \ \       2--torsion\ \ \ \ \ \ \ \ \ \ \ \        3--torsion\ \ \ \ \ \ \ \ \ \      5--torsion\ \ \ \ \ \  7--torsion

\cl{\vrule height .4pt width 4.5in}

$\Gamma_2^1$, $\Gamma_2^2$$\ \ \ \ \ \ \ \ \ \ $(0,6) or (1,2)$\ \ \  \ \ \ \ \ \ \ $(0,4)$\ \ \ \ \ \ \ \ \ \ \ \ \ \ \ \ \ \ $(0,3)$\ \ \ \ \ \ \ \ \ \ $     No

$\Gamma_2^3$ $\ \ \ \ \ \ \ \ \ \ \ \ \ \ $(0,6)$\ \ \ \ \ \ \ \ \ \ \ \ \ \ \ \ \ \ \ \ \ $(0,4)$\ \ \ \ \ \ \ \ \ \ \ \ \ \ \ \ \ $ (0,3)$\ \ \ \ \ \ \ \ \ \ $  No

$\Gamma_2^4$ $\ \ \ \ \ \ \ \ \ \ \ \ \ \ $(0,6)$\ \ \ \ \ \ \ \ \ \ \ \ \ \ \ \ \ \ \ \ \ $(0,4)$\ \ \ \ \ \ \ \ \ \ \ \ \ \ \ \ \ $   No $\ \ \ \ \ \ \ \ \ \  \ \ $  No

$\Gamma_2^5$, $\Gamma_2^6$$\ \ \ \ \ \ \ \ \ \ $(0,6)$\ \ \ \ \ \ \ \ \ \ \ \ \ \ \ \ \ \ \ \ \ $ No $\ \ \ \ \ \ \ \ \ \ \ \ \ \ \ \ \ \ $ No $\ \ \ \ \ \ \ \ \ \ \ \ $No

$\Gamma_3^1$, $\Gamma_3^2$ $\ \ \ \ \ \ \ \ \ $(0,8) or (1,4) $\ \ \ \ \ \ \ \ \ $(1,2) or (0,5)$\ \ \ \ \ \ \ $ No $\ \ \ \ \ \ \ \ \ \ \ \ $(0,3)               

$\Gamma_3^3$$\ \ \ \ \ \ \ \ \ \ \ \ \ \ \ $(0,8) or (1,4)$\ \ \ \ \ \ \ \ \ \ $(0,5)$\ \ \ \ \ \ \ \ \ \ \ \ \ \ \ \ \ $ No$\ \ \ \ \ \ \ \ \ \ \ \ $ (0,3)

$\Gamma_3^4$ $\ \ \ \ \ \ \ \ \ \ \ \ \ \ $(0,8) or (1,4)$\ \ \ \ \ \ \ \ \ \ $(0,5)$\ \ \ \ \ \ \ \ \ \ \ \ \ \ \ \ \ \ \ $No$\ \ \ \ \ \ \ \ \ \ \ $ No

$\Gamma_3^5$ $\ \ \ \ \ \ \ \ \ \ \ \ \ \ $(0,8)$\ \ \ \ \ \ \ \ \ \ \ \ \ \ \ \ \ \ \ \ \ $(0,5)$\ \ \ \ \ \ \ \ \ \ \ \ \ \ \ \ \ \ $No$\ \ \ \ \ \ \ \ \ \ \ \ \ $ No

$\Gamma_3^6$, $\Gamma_3^7$,$\Gamma_3^8$$\ \ \ \ \ \ $(0,8)$\ \ \ \ \ \ \ \ \ \ \ \ \ \ \ \ \ \ \ \ $ No$\ \ \ \ \ \ \ \ \ \ \ \ \ \ \ \ \ \ \ $ No$\ \ \ \ \ \ \ \ \ \ \ \ \ $ No

$\Gamma_1^1$, $\Gamma_1^2$$\ \ \ \ \ \ \ \ \ $ (0,4) $\ \ \ \ \ \ \ \ \ \ \ \ \ \ \ \ \ \ \ $ (0,3) $\ \ \ \ \ \ \ \ \ \ \ \ \ \ \ \ $ No$\ \ \ \ \ \ \ \ \ \ \ \ $ No 

$\Gamma_1^3$ $\ \ \ \ \ \ \ \ \ \  \ \ \ \ $(0,4)$\ \ \ \ \ \ \ \ \ \ \ \ \ \ \ \ \ \ \ \ \ $(0,3)$\ \ \ \ \ \ \ \ \ \ \ \ \ \ \ \ \ $ No$\ \ \ \ \ \ \ \ \ \ \ \ $ No
            
$\Gamma_1^4$ $\ \ \ \ \ \ \ \ \ \ \ \ \ \ $(0,4)$\ \ \ \ \ \ \ \ \ \ \ \ \  \ \ \ \ \ \ \ \  $No $\ \ \ \ \ \ \ \  \ \ \ \ \ \ \ \ \ \ \ $ No$\ \ \ \ \ \ \ \ \ \ \ $ No

\cl{\vrule height .4pt width 4.5in}

}
\end{remark}

\section{\bf{The calculation of the $p$-torsion of the Farrell cohomology of $\Gamma_g^i$ for $g=1,2,3$, $i\ge 1$ and $p$ is an odd prime.}}

Now we begin to analyze $N(\mathbb Z/p)$ in $\Gamma_g^i$. Note that the Riemann Hurwitz equation has two different types of solutions, namely (h,t)=(0,t) or (h,t)=(1,t), where $t$ varies for different $\Gamma_g^i$. We will deal with these two cases separately.

Case 1, $(h,t)=(0,t)$:\qua For $\mathbb Z/p<\Gamma_g^i,$ the $\mathbb Z/p$ action on $S_g$ has $t$ fixed points with quotient space $S_0$. Following arguments similar to those in [4] Lemmas 2.14-2.19, we know that $N(\mathbb Z/p)/\mathbb Z/p$ maps injectively into $\tilde \Gamma_0^t$. In fact, it is not hard to construct this explicit injective mapping. We will give a brief description, but omit the details. 

Any element of $N(\mathbb Z/p)<\Gamma_g^i$ can be lifted to a diffeomorphism of $S_g$. This diffeomorphism has a special property: for the $t$ fixed points of the $\mathbb Z/p$ action, it fixes $i$ of them which associated to the $i$ fixed points in $\Gamma_g^i$; it permutes the other $t-i$. Hence, this diffeomorphism induces a diffeomorphism of the quotient space $S_0$ 
with $t$ points permuted. Thus it gives an element of $\tilde\Gamma_0^t$. The details can be found in [4] or my Ph.D. thesis (1998) at Ohio State University. 

It is well known that $1\to \Gamma_0^t\longrightarrow \tilde \Gamma_0^t\longrightarrow {\sum_t}\to 1,$ where $\sum_t$ is the symmetric group on $t$ letters. Note that in [2], Cohen uses $K_t$ to denote $\Gamma_0^t$. From now on we will adopt his notation $K_t$ for our notation $\Gamma_0^t$. In [2], Cohen calculated $H^*(K_t,\mathbb Z)$ and the action of $\sum_t$ on 
$H^*(K_t,\mathbb Z)$, which implied some cohomology information for $\tilde \Gamma_0^t.$ We will  construct a similar short exact sequence as above which is related to $N(\mathbb Z/p)/\mathbb Z/p$. We can then calculate the Farrell cohomology $\hat H^*(N(\mathbb Z/p)/\mathbb Z/p, F_p)$. 

By the arguments in [4] Lemmas 2.14-2.19, we know that $K_t<N(\mathbb Z/p)/\mathbb Z/p$. Together  with the fact that $N(\mathbb Z/p)/\mathbb Z/p$ maps injectively into $\tilde \Gamma_0^t$, we have a short exact sequence: $1\to K_t\longrightarrow N(\mathbb Z/p)/\mathbb Z/p \longrightarrow {\sum_l}\to 1$, where $\mathbb Z/p< \Gamma_g^i$ and $\sum_l<\sum_t$ is a symmetric group on $l$ letters. The value of $l$ is determined by the fixed point data in the following way: (The details can be found in [4]) Assume that $\alpha\in\Gamma_g^i$ is an element of order $p$. The lifting of $\alpha$ in $\hbox{\it Diffeo}^+(S_g,P_1,...P_i)$ fixes $P_1,P_2,...,P_i,P_{i+1},...P_{t}.$ We denote the fixed point data: $\delta(\alpha)=(1,\beta_2...\beta_i|\beta_{i+1}...\beta_t),$ where $1,\beta_2...\beta_i$ corresponds to $P_1,P_2,...,P_i$ respectively, and $\beta_{i+1}...\beta_t$ corresponds to $P_{i+1},...P_{t}$ respectively. Recall that any element of $N(\mathbb Z/p)<\Gamma_g^i$ can be lifted to an element in $\hbox{\it Diffeo}^+(S_g,P_1,...P_i)$, which is a diffeomorphism fixing $P_1,...P_i$. In [4], we have proved that this diffeomorphism may permute $P_{i+1},...P_{t}$. The value of $l$ is the number of points which are indeed permuted by the diffeomorphism, it is determined by the fixed point data of $\alpha$.(Note that $l$ is at most $t-i$.) We will use some examples to illustrate how to determine $l$.  E.g., 1) $\delta(\alpha)=(1,1|2,1,1).$ Here $i=2$ and $t=5$. Any element of $N(\mathbb Z/p)$ (We abuse the notation an element with its lifting diffeomorphism) fixes $P_1, P_2$ and may permute $P_3,P_4,P_5.$ By Lemma 2.16 in [4], if the $P_j$ and the $P_k$ are allowed to be permuted by any  element of $N(\mathbb Z/p)$, then $\beta_j=\beta_k$. Therefore, $l\le 2$. Indeed $l=2$ (See details in [4]). E.g., 2) $\delta(\alpha)=(1,2|1,1,1).$ Then $l=3$. (The proof can be found in [4]) Knowing the value of $l$, the above short exact sequence is completely determined. In Cohen's paper [2], we can find explicitly the $\sum_t$ action on $H^*(K_t,\mathbb Z)$. Hence the $\sum_l$ action on $H^*(K_t,\mathbb Z)$ is known. Now we can apply the Serre spectral sequence with respect to the above short exact sequence to calculate the Farrell cohomology $\hat H^*(N(\mathbb Z/p)/\mathbb Z/p,F_p)$. 

In order to get $\hat H^*(N(\mathbb Z/p),\mathbb Z)_{(p)}$, we need to consider another short exact sequence $1\to \mathbb Z/p\to N(\mathbb Z/p)\to N(\mathbb Z/p)/\mathbb Z/p\to 1.$ This short exact sequence is central. Thus the associated Serre spectral sequence has trivial coefficients, so we can calculate $\hat H^*(N(\mathbb Z/p),\mathbb Z)_{(p)}$. In fact, the above central property comes from the periodicity of pure mapping class groups with punctures. In [4], we proved that $N(\mathbb Z/p)=C(\mathbb Z/p)$ for any $\mathbb Z/p<\Gamma_g^i$. It is a corollary of periodicity. Note that it is in contrast to the unpunctured mapping class groups, which are not periodic in general.  

Case 2, $(h,t)=(1,t)$:\qua For $\mathbb Z/p<\Gamma_g^i,$ the $\mathbb Z/p$ action on $S_g$ has $t$ fixed points and the quotient space is $S_1$. By [5], $N(\mathbb Z/p)/\mathbb Z/p$ 
can be viewed as a subgroup of the mapping class group $\tilde \Gamma_1^t$ of finite index. Xia  in [7] developed a way to calculate $\hat H^*(N(\mathbb Z/p),\mathbb Z)_{(p)}$ in this case, which we can adapt for our cases. However, in his case, the period of his mapping class group is 4, whereas in our case the period is always $2$. We will show later in this paper where his result does not apply to our cases, and which modifications are necessary. 

Now in both of the above two cases, we can calculate $\hat H^*(N(\mathbb Z/p),\mathbb Z)_{(p)}.$ In order  to apply Brown's theorem mentioned in the introduction, we need to count the conjugacy classes of subgroup of order $p$ in $\Gamma_g^i$. The tools we use are Theorem 2.5 and Proposition 2.6 in [4]:

\begin{theorem}[Theorem 2.5 in \cite{4}]
Let $\Gamma^i_g =\pi _0{(\hbox{\it Diffeo}^+(S_g,P_1,...,P_i))}$, and let $\alpha,\alpha'$ be elements of order $p$ in $\Gamma^i_g$, with $\delta_i (\alpha)=
( \beta_1,...,\beta_i|\beta_{i+1},...,\beta_t)$, 
 $\delta_i (\alpha')=( \beta_1',...,\beta_i'|\beta_{i+1}',...\beta_t')$. Then, the following holds:

The element $\alpha$ is conjugate to $\alpha'$ in $\Gamma^i_g$ if and only if $\beta_1=\beta_1'$, ..., $\beta_i=\beta_i'$, and $(\beta_{i+1},...,\beta_t)=(\beta_{i+1}',...
,\beta_t')$ as unordered integer tuples; i.e.,
two elements of order $p$ in $\Gamma^i_g$ are conjugate if and only if they have
the same fixed point data.
\end{theorem}   

\begin{proposition}[Proposition 2.6 in \cite{4}]
Let $t$ be a non-negative integer\break which satisfies the Riemann Hurwitz equation 
$2g-2=p(2h-2)+t(p-1)$ with $t\neq 1$ and $t\ge i$. Then the number of different integer
 tuples $(1,\beta_2,...,\beta_i|$ $\beta_{i+1},...,\beta_t)$ such that $(1,\beta_2,...,\beta_i)$ is ordered, $(\beta_{i+1},...,\beta_t)$ is unordered, and  
$1+\beta_2+...+\beta_t=0(\mod p)$, where $0<\beta_i<p$ for all $i$, is the 
same as the number of conjugacy classes of subgroups of order $p$ in 
$\Gamma^i_g$ which act on $S_g$ with $t$ fixed points. 
\end{proposition}

Here we know that the number of conjugacy classes is the number of different integer t-tuples $(1,\beta_2,...,\beta_i|\beta_{i+1},...\beta_t)$ such that $(1,\beta_2,...\beta_i)$ is ordered, $(\beta_{i+1},...\beta_t)$ is unordered, and $1+\beta_2+...+\beta_t=0(\mod p)$, where $0<\beta_j<p$ for all $j$. Each solution of the above equation corresponds to one type of fixed point data, which then corresponds to one conjugacy class of subgroups of order $p$. To get the integer t-tuple is a simple algebraic problem which we will not cover. Instead, we will give the fixed point data types directly in the following theorems without details. 

In general, the above arguments work for any prime $p$. However, the $2-$primary component of Farrell cohomology is very complicated for the cases $\Gamma_2^i$ and $\Gamma_3^i$. Its calculation remains an open question.

\begin{theorem}
\label{e}
\begin{enumerate}\romanIndex 
\item $\hat H^i(\Gamma^1_1,\mathbb Z)$ is known, because $$\Gamma_1^1=SL_2(\mathbb Z)\cong \mathbb Z/4*_{\mathbb Z/2}\mathbb Z/6.$$
\item 
$$
\hat H^i(\Gamma_1^2,\mathbb Z)_{(2)}=\left\{
\aligned \mathbb Z/4  \ \ \ \ \ \ \ \ i=0~mod(2)\\
\mathbb Z/2  \ \ \ \ \ \ \ \ i=1~mod(2)
\endaligned
\right\}.
$$
$$
\hat H^i(\Gamma_1^2,\mathbb Z)_{(3)}=\left\{
\aligned \mathbb Z/3  \ \ \ \ \ \ \ \ i=0~mod(2)\\
0  \ \ \ \ \ \ \ \ i=1~mod(2)
\endaligned
\right\}.
$$
\item
$$
\hat H^i(\Gamma_1^3,\mathbb Z)_{(2)}=\left\{
\aligned \mathbb Z/2  \ \ \ \ \ \ \ \ i=0~mod(2)\\
\mathbb Z/2\oplus \mathbb Z/2  \ \ \ \ \ \ \ \ i=1~mod(2)
\endaligned
\right\}.
$$
$$
\hat H^i(\Gamma_1^3,\mathbb Z)_{(3)}=\left\{
\aligned \mathbb Z/3  \ \ \ \ \ \ \ \ i=0~mod(2)\\
0  \ \ \ \ \ \ \ \ i=1~mod(2)
\endaligned
\right\}.
$$
\item
$$
\hat H^i(\Gamma_1^4,\mathbb Z)_{(2)}=\left\{
\aligned \mathbb Z/2  \ \ \ \ \ \ \ \ i=0~mod(2)\\
\mathbb Z/2\oplus \mathbb Z/2  \ \ \ \ \ \ \ \ i=1~mod(2)
\endaligned
\right\}.
$$
\item
$$
\hat H^i(\Gamma_1^i,\mathbb Z)=0, for\ \ i\ge 5.
$$
\end{enumerate}
\end{theorem}

\proof
Case (1): $p=3$

Since $p=3$ and $g=(p-1)/2=1$, by [4] Theorem 2.23,
 
$$
\hat H^i(\Gamma_1^2,\mathbb Z)_{(3)}=\left\{
\aligned \mathbb Z/3  \ \ \ \ \ \ \ \ i=0~mod(2)\\
0  \ \ \ \ \ \ \ \ i=1~mod(2)
\endaligned
\right\}.
$$
$$
\hat H^i(\Gamma_1^3,\mathbb Z)_{(3)}=\left\{
\aligned \mathbb Z/3  \ \ \ \ \ \ \ \ i=0~mod(2)\\
0  \ \ \ \ \ \ \ \ i=1~mod(2)
\endaligned
\right\}.
$$

Case (2): $p=2$

From Corollary \ref{d}, we know that for $\Gamma_1^2$, $\Gamma_1^3$, $\Gamma_1^4$, $(h,t)=(0,4)$. Thus the $\mathbb Z/2$ action on $S_1$ has $4$ fixed points and the quotient space is $S_0$. 
We have the following short exact sequences:

(i)\qua For $\Gamma_1^3$ or $\Gamma_1^4$, $1\to K_4\longrightarrow N(\mathbb Z/2)/\mathbb
Z/2 \longrightarrow {\sum_1}\to 1$, where $\mathbb Z/2< \Gamma_1^3$ or $\Gamma_1^4$. The  corresponding fixed point data is $(1,1,1|1)$ or $(1,1,1,1|)$ respectively.

(ii)\qua For $\Gamma_1^2$, $1\to K_4\longrightarrow N(\mathbb Z/2)/\mathbb
Z/2 \longrightarrow {\sum_2}\to 1$, where $\mathbb Z/2< \Gamma_1^2.$ The corresponding fixed point data is $(1,1|1,1)$.

Case (2)(i):\qua For $\Gamma_1^3$, $1\to K_4\longrightarrow N(\mathbb Z/2)/\mathbb Z/2 \longrightarrow {\sum_1}\to 1$. It is known (Cohen [2]), that  
$$ H^i(N(\mathbb Z/2)/\mathbb Z/2,\mathbb Z)\cong H^i(K_4,\mathbb Z)\cong \left\{
\aligned \mathbb Z  \ \ \ \ \ \ \ \ i=0\\
\mathbb Z\oplus\mathbb Z  \ \ \ \ \ \ \ \ i=1\\
0  \ \ \ \ \ \ \ \ i\ge 2
\endaligned
\right\}.
$$
It is easy to see that $H^2(N(\mathbb Z/2)/\mathbb Z/2,F_2)=0.$ So the following short exact sequence splits: $1\to \mathbb Z/2\longrightarrow N(\mathbb Z/2) \longrightarrow N(\mathbb Z/2)/\mathbb Z/2\to 1$. Also, this short exact sequence is central. In fact, $N(\mathbb Z/p)=C(\mathbb Z/p)$ for any $\mathbb Z/p<\Gamma_g^i$. (It is because of the periodicity we proved in [4])

We have $N(\mathbb Z/2)\cong N(\mathbb Z/2)/(\mathbb Z/2)\times \mathbb Z/2$.

\noindent By the K\"unneth theorem,
 
$\hat H^i(N(\mathbb Z/2),\mathbb Z)_{(2)}=\mathbb Z/2  \ \ \ \ \ \ \ \ i=0~mod(2)$

$\hat H^i(N(\mathbb Z/2),\mathbb Z)_{(2)}=\mathbb Z/2\oplus \mathbb Z/2  \ \ \ \ \ \ \ \ i=1~mod(2)$

\noindent Since there is one type of fixed point data, namely $(1,1,1|1)$, we have 

$\hat H^i(\Gamma_1^3,\mathbb Z)_{(2)}=\mathbb Z/2  \ \ \ \ \ \ \ \ i=0~mod(2)$

$\hat H^i(\Gamma_1^3,\mathbb Z)_{(2)}=\mathbb Z/2\oplus \mathbb Z/2  \ \ \ \ \ \ \ \ i=1~mod(2)$.

It is similar for the case of $\Gamma_1^4$.

Case (2)(ii):\qua For $\Gamma_1^2$, $1\to K_4\longrightarrow N(\mathbb Z/2)/\mathbb Z/2 \longrightarrow {\sum_2}\to 1$.

The Serre spectral sequence takes the form 
 
$E^{ij}_2=H^i(\Sigma_2,H^j(K_4,F_2))\Longrightarrow H^{i+j}(N(\mathbb Z/2)/\mathbb Z/2,F_2)$. Recall that in case(2), we always have $(h,t)=(0,4)$, any $\mathbb Z/2$(the lifting of $\mathbb Z/2$) acts on $S_1$ with 4 fixed points. Let us assume that the four 
fixed points are $P_1,P_2,P_3,P_4.$ Since the fixed point data is $(1,1|1,1)$, we know that the elements of the normalizer $N(\mathbb Z/2)$ fix two points and permute the other two. Without loss of generality, $P_1,P_2$ are fixed and $P_3, P_4$ are permuted by the elements of the normalizer $N(\mathbb Z/2)$. So, $\sum_2$ 
is generated by $<x>\ \ =\ \ <(34)>$. We need to calculate $H^i(<x>, H^j(K_4,F_2))$. This is related to $H^j(K_4,F_2)_{<x>}$ and $H^j(K_4,F_2)^{<x>}$, the coinvariant and invariant of $x$ on $H^*(K_4,F_2)$. By Cohen in [2], $H^1(K_4,F_2)$ is generated by two degree-one generators $\{B_{42}, B_{43}\}$. Also, $(34)B_{42}=-B_{42}$, $(34)B_{43}=B_{42}+B_{43}$. So the invariant is generated by $B_{42}$.
Thus we have 

$H^1(K_4,F_2)^{<x>}=<B_{42}>\cong F_2.$

Similarly, $H^1(K_4,F_2)_{<x>}=H^1(K_4,F_2)/M$, where $M=<y-xy>$ and $y\in H^1(K_4,F_2)$.
So, $$H^1(K_4,F_2)_{<x>}=<\bar B_{43}>\cong F_2.$$

Consider the norm map $N:H^1(K_4,F_2)_{<x>}\rightarrow H^1(K_4,F_2)^{<x>}$. It is easy to verify that $N(\bar B_{43})=B_{42}$, 
Hence, $N$ is an isomorphism. Therefore,
$$H^i(<x>,H^1(K_4,F_2))=coker\ \ N=0, \ \ if\ \ i=0~mod(2)(i>0),$$
$$H^i(<x>,H^1(K_4,F_2))=ker\ \ N=0, \ \ if\ \ i=1~mod(2),$$
$$H^0(<x>,H^1(K_4,F_2))=H^1(K_4,F_2)^{<x>}\cong F_2.$$

Also, since $H^j(K_4,F_2)=0 \ \ for\ \ j\ge 2,$ we have $$H^i(<x>, H^j(K_4,F_2))=0, \ \ for\ \ j\ge 2,$$
$$H^i(<x>, H^0(K_4,F_2))\cong H^i(<x>, F_2)\cong F_2.$$

Since the Serre spectral sequence collapses, we have $H^{0}(N(\mathbb Z/2)/\mathbb Z/2,F_2)\cong F_2,$ $H^{1}(N(\mathbb Z/2)/\mathbb Z/2,F_2)\cong F_2\oplus F_2,$ and $H^{i}(N(\mathbb Z/2)/\mathbb Z/2,F_2)\cong F_2, for\ \ i\ge 2$

Now we need to analyze the short exact sequence $1\to \mathbb Z/2\longrightarrow N(\mathbb Z/2) \longrightarrow N(\mathbb Z/2)/\mathbb Z/2\to 1$. We claim that in this case there is $\mathbb Z/4$ in $N(\mathbb Z/2)$. The reason is the following.

First, we show that there is $\mathbb Z/4$ in $\Gamma_1^2$. In order to detect if there is a $\mathbb Z/p^2$ in $\Gamma_1^2$, we need to use the generalized Riemann Hurwitz equation(pg.259 [3]): $2g-2=p^2(2h-2)+sp^2(1-1/p)+tp^2(1-1/p^2)$, where $s$ is the number of order $p$ singular points and $t$ is the number of order $p^2$ singular points. Here we have $g=1$ and $p=2$, the equation has a solution $(h,s,t)=(0,1,2).$ By similar reasons we know that there is a $\mathbb Z/4$ in $\Gamma_1^2$. In fact, $\mathbb Z/4$(the lifting of $\mathbb Z/4$) acts on $S_1$ with two singular points of order 4( because $t=2$), one singular point of order 2 (because $s=1$) and the orbit space is $S_1/(\mathbb Z/4)=S_0.$(because $h=0$). Note that there is no $\mathbb Z/4$ in $\Gamma_1^3$ or $\Gamma_1^4$. This is because any lifting of $\mathbb Z/4$ acting on $S_1$ must fix at least three points,  which contradicts two singular points of order 4. (A fixed point is a special case of a singular point, where the stablizer is the entire group that acts.)

Since there is only one kind of fixed point data corresponding to $\mathbb Z/2<\Gamma_1^2$, namely $(1,1|1,1)$, the above $\mathbb Z/4$ must be contained in $N(\mathbb Z/2)$.

Now look at the Serre spectral sequence associated to the short exact sequence: 
 $$1\to \mathbb Z/2\longrightarrow N(\mathbb Z/2) \longrightarrow N(\mathbb Z/2)/\mathbb Z/2\to 1.$$ 
Compare it with the Serre spectral sequence associated to the short exact sequence:
$$1\to \mathbb Z/2\longrightarrow \mathbb Z/4\longrightarrow \mathbb Z/2\to 1.$$
In fact, $N(\mathbb Z/2)$ is periodic,(by the result that $\Gamma_1^2$ is periodic in [4]), so there is no $\mathbb Z/2\times \mathbb Z/2$ in $\Gamma_1^2$. We have,  

$\hat H^i(N(\mathbb Z/2),\mathbb Z)_{(2)}=\mathbb Z/4  \ \ \ \ \ \ \ \ i=0~mod(2),$

$\hat H^i(N(\mathbb Z/2),\mathbb Z)_{(2)}=\mathbb Z/2  \ \ \ \ \ \ \ \ i=1~mod(2).$

There is only one conjugacy class of $\mathbb Z/2$ in $\Gamma_1^2$(one kind of fixed point data), so we have 
$$
\hat H^i(\Gamma_1^2,\mathbb Z)_{(2)}=\left\{
\aligned \mathbb Z/4  \ \ \ \ \ \ \ \ i=0~mod(2)\\
\mathbb Z/2  \ \ \ \ \ \ \ \ i=1~mod(2)
\endaligned
\right\}.\eqno{\qed}
$$

In the remainder of this paper we let $n\mathbb Z/p$ denote the direct sum of $\mathbb Z/p$ with itself $n$ times. 

\begin{theorem}
\label{f}
\begin{enumerate}\romanIndex 
\item 
$$
\hat H^i(\Gamma_2^1,\mathbb Z)_{(3)}=\left\{
\aligned \mathbb Z/3\oplus\mathbb Z/3  \ \ \ \ \ \ \ \ i=0~mod(2)\\
\mathbb Z/3  \ \ \ \ \ \ \ \ i=1~mod(2)
\endaligned
\right\}.
$$
$$
\hat H^i(\Gamma_2^1,\mathbb Z)_{(5)}=\left\{
\aligned \mathbb Z/5\oplus\mathbb Z/5  \ \ \ \ \ \ \ \ i=0~mod(2)\\
0  \ \ \ \ \ \ \ \ i=1~mod(2)
\endaligned
\right\}.
$$
\item
$$
\hat H^i(\Gamma_2^2,\mathbb Z)_{(3)}=\left\{
\aligned 3\mathbb Z/3  \ \ \ \ \ \ \ \ i=0~mod(2)\\
3\mathbb Z/3  \ \ \ \ \ \ \ \ i=1~mod(2)
\endaligned
\right\}.
$$
$$
\hat H^i(\Gamma_2^2,\mathbb Z)_{(5)}=\left\{
\aligned 3\mathbb Z/5  \ \ \ \ \ \ \ \ i=0~mod(2)\\
0  \ \ \ \ \ \ \ \ i=1~mod(2)
\endaligned
\right\}.
$$
\item
$$
\hat H^i(\Gamma_2^3,\mathbb Z)_{(3)}=\left\{
\aligned 3\mathbb Z/3  \ \ \ \ \ \ \ \ i=0~mod(2)\\
6\mathbb Z/3  \ \ \ \ \ \ \ \ i=1~mod(2)
\endaligned
\right\}.
$$
$$
\hat H^i(\Gamma_2^3,\mathbb Z)_{(5)}=\left\{
\aligned 3\mathbb Z/5  \ \ \ \ \ \ \ \ i=0~mod(2)\\
0  \ \ \ \ \ \ \ \ i=1~mod(2)
\endaligned
\right\}.
$$
\item
$$
\hat H^i(\Gamma_2^4,\mathbb Z)_{(3)}=\left\{
\aligned 3\mathbb Z/3  \ \ \ \ \ \ \ \ i=0~mod(2)\\
6\mathbb Z/3  \ \ \ \ \ \ \ \ i=1~mod(2)
\endaligned
\right\}.
$$
\item
$$
\hat H^*(\Gamma_2^i,\mathbb Z)=0, for\ \ i\ge 7.
$$
\end{enumerate}
\end{theorem}
\begin{proof}

Case (1): $p=5$

Since $p=5$ and $g=(p-1)/2=2$, by [4] Theorem 2.23,
 
$$
\hat H^i(\Gamma_2^1,\mathbb Z)_{(5)}=\left\{
\aligned 2 \mathbb Z/5  \ \ \ \ \ \ \ \ i=0~mod(2)\\
0  \ \ \ \ \ \ \ \ i=1~mod(2)
\endaligned
\right\}.
$$
$$
\hat H^i(\Gamma_2^2,\mathbb Z)_{(5)}=\left\{
\aligned 3 \mathbb Z/5  \ \ \ \ \ \ \ \ i=0~mod(2)\\
0  \ \ \ \ \ \ \ \ i=1~mod(2)
\endaligned
\right\}.
$$
$$
\hat H^i(\Gamma_2^3,\mathbb Z)_{(5)}=\left\{
\aligned 3 \mathbb Z/5  \ \ \ \ \ \ \ \ i=0~mod(2)\\
0  \ \ \ \ \ \ \ \ i=1~mod(2)
\endaligned
\right\}.
$$

Case (2): $p=3$

From Corollary \ref{b}, we know that for $\Gamma_2^1$, $\Gamma_2^2$, $\Gamma_2^3$, $\Gamma_2^4$, $(h,t)=(0,4)$. The $\mathbb Z/3$ action on $S_2$ has $4$ fixed points.
Following arguments similar to those in [4] Lemmas 2.14-2.19, we have short exact sequences:

(i)\qua For $\Gamma_2^3$ or $\Gamma_2^4$, $1\to K_4\longrightarrow N(\mathbb Z/3)/\mathbb
Z/3 \longrightarrow {\sum_1}\to 1$, where $\mathbb Z/3< \Gamma_2^3\ \ or\ \ \Gamma_2^4;$ the corresponding fixed point data is $(1,1,2|2)$, or $(1,2,1|2)$, or $(1,2,2|1)$ for $\mathbb Z/3< \Gamma_2^3$, and $(1,1,2,2|)$, or $(1,2,1,2|)$, or $(1,2,2,1|)$ for $\mathbb Z/3< \Gamma_2^4$.

(ii)(a)\qua For $\Gamma_2^2$, $1\to K_4\longrightarrow N(\mathbb Z/3)/\mathbb
Z/3 \longrightarrow {\sum_2}\to 1$, where $\mathbb Z/3< \Gamma_2^2$ and the corresponding fixed point data for $\mathbb Z/3$ is $(1,1|2,2).$

(ii)(b)\qua For $\Gamma_2^2$, $1\to K_4\longrightarrow N(\mathbb Z/3)/\mathbb
Z/3 \longrightarrow {\sum_1}\to 1$, where $\mathbb Z/3< \Gamma_2^2$ and the corresponding fixed point data for $\mathbb Z/3$ is $(1,2|1,2).$

(iii)\qua For $\Gamma_2^1$, $1\to K_4\longrightarrow N(\mathbb Z/3)/\mathbb
Z/3 \longrightarrow {\sum_2}\to 1$, where $\mathbb Z/3< \Gamma_2^1$ and the fixed point data for $\mathbb Z/3$ is $(1|1,2,2).$

Case(2)(i):\qua For $\Gamma_2^3$, $1\to K_4\longrightarrow N(\mathbb Z/3)/\mathbb Z/3 \longrightarrow {\sum_1}\to 1$. 

As in case (2)(i) in Theorem \ref{e},

$\hat H^i(N(\mathbb Z/3),\mathbb Z)_{(3)}=\mathbb Z/3  \ \ \ \ \ \ \ \ i=0~mod(2)$

$\hat H^i(N(\mathbb Z/3),\mathbb Z)_{(3)}=\mathbb Z/3\oplus \mathbb Z/3  \ \ \ \ \ \ \ \ i=1~mod(2)$

We have three different types of the fixed point data, namely $(1,1,2|2)$, 

$(1,2,1|2)$, and $(1,2,2|1)$. Therefore, 

$\hat H^i(\Gamma_2^3,\mathbb Z)_{(3)}=3\mathbb Z/3  \ \ \ \ \ \ \ \ i=0~mod(2)$

$\hat H^i(\Gamma_2^3,\mathbb Z)_{(3)}=6\mathbb Z/3  \ \ \ \ \ \ \ \ i=1~mod(2)$.

It is similar for the case of $\Gamma_2^4$.

Case(2)(ii)(a):\qua For $\Gamma_2^2$, if the fixed point data is $(1,1|2,2)$, there is a short exact sequence: $1\to K_4\longrightarrow N(\mathbb Z/3)/\mathbb Z/3 \longrightarrow {\sum_2}\to 1$.

The Serre spectral sequence takes the form 

$E^{ij}_2\cong H^i(\Sigma_2,H^j(K_4,F_3))\Longrightarrow H^{i+j}(N(\mathbb Z/3)/\mathbb Z/3,F_3)$. 

Since 2 and 3 are relatively prime, $E^{ij}_2=0$ for $i>0$. Thus we only need to consider $i=0.$ 
$E^{0j}_2=H^j(K_4,F_3)^{\Sigma_2}.$ 
Let us assume that the four 
fixed points of the $\mathbb Z/3$ action on $S_2$ are $P_1,P_2,P_3,P_4.$ Since the fixed point data is $(1,1|2,2)$, we know that the elements of the normalizer $N(\mathbb Z/3)$ fix two points and permute the other two. Without loss of generality, $P_1,P_2$ are fixed and $P_3, P_4$ are permuted by the elements of the normalizer $N(\mathbb Z/3)$ . So $\sum_2$ 
is generated by $<x>\ \ =\ \ <(34)>$. By Cohen in [2], $H^1(K_4,F_3)$ is generated by two degree-one generators $\{B_{42}, B_{43}\}$. As in the proof of case(2)(ii) in Theorem \ref{e}, we get $$H^1(K_4,F_3)^{<x>}=<B_{42}+2B_{43}>\cong F_3.$$
$$H^j(K_4,F_3)^{<x>}=0\ \ for\ \ j\ge 2.$$
$$H^0(K_4,F_3)^{<x>}\cong F_3.$$

By the Serre spectral sequence, we have $H^{0}(N(\mathbb Z/3)/\mathbb Z/3,F_3)\cong F_3,$ 

$H^{1}(N(\mathbb Z/3)/\mathbb Z/3,F_3)\cong F_3,$ 

and $H^{i}(N(\mathbb Z/3)/\mathbb Z/3,F_3)=0\ \ for\ \ i>1.$

Now look at the spectral sequence associated to the short exact sequence: 
 $$1\to \mathbb Z/3\longrightarrow N(\mathbb Z/3) \longrightarrow N(\mathbb Z/3)/\mathbb Z/3\to 1.$$ 

As in Case (2)(i) in Theorem 2.3, using the K\"unneth theorem,
 
$\hat H^i(N(\mathbb Z/3),\mathbb Z)_{(3)}=\mathbb Z/3\oplus\mathbb Z/3  \ \ \ \ \ \ \ \ i=0~mod(2)$

$\hat H^i(N(\mathbb Z/3),\mathbb Z)_{(3)}=\mathbb Z/3  \ \ \ \ \ \ \ \ i=1~mod(2)$

Case (2)(ii)(b):\qua For the fixed point data $(1,2|1,2)$, the short exact sequence is  
$1\to K_4\longrightarrow N(\mathbb Z/3)/\mathbb Z/3 \longrightarrow {\sum_1}\to 1.$ Hence, as in  case (2)(i), we have 

$\hat H^i(N(\mathbb Z/3),\mathbb Z)_{(3)}=\mathbb Z/3  \ \ \ \ \ \ \ \ i=0~mod(2)$

$\hat H^i(N(\mathbb Z/3),\mathbb Z)_{(3)}=\mathbb Z/3\oplus \mathbb Z/3  \ \ \ \ \ \ \ \ i=1~mod(2)$

Now put case (2)(ii)(a) and case (2)(ii)(b) together to get 

$\hat H^i(\Gamma_2^2,\mathbb Z)_{(3)}=3\mathbb Z/3  \ \ \ \ \ \ \ \ i=0~mod(2)$

$\hat H^i(\Gamma_2^2,\mathbb Z)_{(3)}=3\mathbb Z/3  \ \ \ \ \ \ \ \ i=1~mod(2)$.

Case (2)(iii):\qua As in case (ii)(a), we have 

$\hat H^i(N(\mathbb Z/3),\mathbb Z)_{(3)}=\mathbb Z/3\oplus\mathbb Z/3  \ \ \ \ \ \ \ \ i=0~mod(2)$

$\hat H^i(N(\mathbb Z/3),\mathbb Z)_{(3)}=\mathbb Z/3  \ \ \ \ \ \ \ \ i=1~mod(2)$

Since there is only one conjugacy classes of subgroups of order $p$ corresponding to $(1|1,2,2)$, we have 

$\hat H^i(\Gamma_2^1,\mathbb Z)_{(3)}=2 \mathbb Z/3  \ \ \ \ \ \ \ \ i=0~mod(2)$

$\hat H^i(\Gamma_2^1,\mathbb Z)_{(3)}=\mathbb Z/3  \ \ \ \ \ \ \ \ i=1~mod(2)$.

\end{proof}

\begin{theorem}
\label{g}
\begin{enumerate}\romanIndex 
\item 

$$
\hat H^i(\Gamma_3^1,\mathbb Z)_{(3)}=\left\{
\aligned  (3\mathbb Z/3)\oplus(2\mathbb Z/9) \ \ \ \ \ \ \ \ i=0~mod(2)\\
 2\mathbb Z/3 \ \ \ \ \ \ \ \ i=1~mod(2)
\endaligned
\right\}.
$$
$$
\hat H^i(\Gamma_3^1,\mathbb Z)_{(7)}=\left\{
\aligned 3\mathbb Z/7  \ \ \ \ \ \ \ \ i=0~mod(2)\\
0  \ \ \ \ \ \ \ \ i=1~mod(2)
\endaligned
\right\}.
$$
\item
$$
\hat H^i(\Gamma_3^2,\mathbb Z)_{(3)}=\left\{
\aligned  (6\mathbb Z/3)\oplus\mathbb Z/9\ \ or\ \ (4\mathbb Z/3)\oplus(2\mathbb Z/9)  \ \ \ \ \ \ i=0~mod(2)\\
4\mathbb Z/3  \ \ \ \ \ \ \ \ i=1~mod(2)
\endaligned
\right\}.
$$
$$
\hat H^i(\Gamma_3^2,\mathbb Z)_{(7)}=\left\{
\aligned 5\mathbb Z/7  \ \ \ \ \ \ \ \ i=0~mod(2)\\
0  \ \ \ \ \ \ \ \ i=1~mod(2)
\endaligned
\right\}.
$$
\item
$$
\hat H^i(\Gamma_3^3,\mathbb Z)_{(3)}=\left\{
\aligned  19\mathbb Z/3 \ \ or \ \ (17\mathbb Z/3)\oplus\mathbb Z/9\\
 \ \ or\ \ (15\mathbb Z/3)\oplus(2\mathbb Z/9)\ \ or \ \ (13\mathbb Z/3)\oplus(3\mathbb Z/9)\\
 \ \ or\ \ (11\mathbb Z/3)\oplus(4\mathbb Z/9)  \ \ \ \ \ \ i=0~mod(2)\\
 14\mathbb Z/3   \ \ \ \ \ \ \ \ i=1~mod(2)
\endaligned
\right\}.
$$

$$
\hat H^i(\Gamma_3^3,\mathbb Z)_{(7)}=\left\{
\aligned 5\mathbb Z/7  \ \ \ \ \ \ \ \ i=0~mod(2)\\
0  \ \ \ \ \ \ \ \ i=1~mod(2)
\endaligned
\right\}.
$$

\item
$$
\hat H^i(\Gamma_3^4,\mathbb Z)_{(3)}=\left\{
\aligned   35\mathbb Z/3, \ \ or \ \ (33\mathbb Z/3)\oplus\mathbb Z/9,\ \ or\ \ (31\mathbb Z/3)\oplus(2\mathbb Z/9),\\
\ \ or \ \ (29\mathbb Z/3)\oplus(3\mathbb Z/9), \ \ or \ \ (27\mathbb Z/3)\oplus(4\mathbb Z/9),\\   \ \ or \ \ (25\mathbb Z/3)\oplus(5\mathbb Z/9),\ \ \ \ \ \ \ \ i=0~mod(2)\\
 25\mathbb Z/3   \ \ \ \ \ \ \ \ i=1~mod(2)
\endaligned
\right\}.
$$
\item
$$
\hat H^i(\Gamma_3^5,\mathbb Z)_{(3)}=\left\{
\aligned   35\mathbb Z/3, \ \ or \ \ (33\mathbb Z/3)\oplus\mathbb Z/9,\ \ or\ \ (31\mathbb Z/3)\oplus(2\mathbb Z/9),\\
\ \ or \ \ (29\mathbb Z/3)\oplus(3\mathbb Z/9), \ \ or \ \ (27\mathbb Z/3)\oplus(4\mathbb Z/9),\\
   \ \ or \ \ (25\mathbb Z/3)\oplus(5\mathbb Z/9),\ \ \ \ \ \ \ \ i=0~mod(2)\\
 25\mathbb Z/3   \ \ \ \ \ \ \ \ i=1~mod(2)
\endaligned
\right\}.
$$
\item
$$
\hat H^*(\Gamma_3^i,\mathbb Z)=0, for\ \ i\ge 9.
$$
\end{enumerate}
\end{theorem}
\begin{proof}

Case (1): $p=7$   

Since $p=7$ and $g=(p-1)/2=3$, by [4] Theorem 2.23,
 
$$
\hat H^i(\Gamma_3^1,\mathbb Z)_{(7)}=\left\{
\aligned 3 \mathbb Z/7  \ \ \ \ \ \ \ \ i=0~mod(2)\\
0  \ \ \ \ \ \ \ \ i=1~mod(2)
\endaligned
\right\}.
$$
$$
\hat H^i(\Gamma_3^2,\mathbb Z)_{(7)}=\left\{
\aligned 5 \mathbb Z/7  \ \ \ \ \ \ \ \ i=0~mod(2)\\
0  \ \ \ \ \ \ \ \ i=1~mod(2)
\endaligned
\right\}.
$$
$$
\hat H^i(\Gamma_3^3,\mathbb Z)_{(7)}=\left\{
\aligned 5 \mathbb Z/7  \ \ \ \ \ \ \ \ i=0~mod(2)\\
0  \ \ \ \ \ \ \ \ i=1~mod(2)
\endaligned
\right\}.
$$

Case (2): $p=3$

From Corollary \ref{c}, we know that for $\Gamma_3^3$, $\Gamma_3^4$, $\Gamma_3^5$, $(h,t)=(0,5)$. The $\mathbb Z/3$ action on $S_3$ has $5$ fixed points. 
For $\Gamma_3^1$, $\Gamma_3^2,$ $(h,t)=(0,5)$, or $(h,t)=(1,2).$ The $\mathbb Z/3$ action on $S_3$ has $5$ fixed points with quotient space $S_0$ or the $\mathbb Z/3$ action on $S_3$ has $2$ fixed points with quotient space $S_1.$ This depends on $\mathbb Z/3$'s fixed point data.

As in [4] Lemmas 2.14-2.19, we have short exact sequences:

(i)\qua For $\Gamma_3^4$ or $\Gamma_3^5$, $1\to K_5\longrightarrow N(\mathbb Z/3)/\mathbb
Z/3 \longrightarrow {\sum_1}\to 1$, where $\mathbb Z/3< \Gamma_3^4\ \ or\ \ \Gamma_3^5$ and the fixed point data is $(1,2,2,2|2)$, or $(1,2,1,1|1)$, or $(1,1,2,1|1)$, or $(1,1,1,2|1)$, or $(1,1,1,1|2)$ for $\mathbb Z/3< \Gamma_3^4$, and $(1,2,2,2,2|)$, or $(1,2,1,1,1|)$, or $(1,1,2,1,1|)$, or $(1,1,1,2,1|)$, or $(1,1,1,1,2|)$ for $\mathbb Z/3< \Gamma_3^5$. 

(ii)(a)\qua For $\Gamma_3^3$, $1\to K_5\longrightarrow N(\mathbb Z/3)/\mathbb
Z/3 \longrightarrow {\sum_2}\to 1$, where $\mathbb Z/3< \Gamma_3^3$ and the fixed point data for $\mathbb Z/3$ is $(1,2,2|2,2)$, or $(1,2,1|1,1)$, or $(1,1,2|1,1)$.

(ii)(b)\qua For $\Gamma_3^3$, $1\to K_5\longrightarrow N(\mathbb Z/3)/\mathbb
Z/3 \longrightarrow {\sum_1}\to 1$, where $\mathbb Z/3< \Gamma_3^3$ and the fixed point data for $\mathbb Z/3$ is $(1,1,1|1,2)$.

(iii)(a)\qua For $\Gamma_3^2$, $1\to K_5\longrightarrow N(\mathbb Z/3)/\mathbb
Z/3 \longrightarrow {\sum_3}\to 1$, where $\mathbb Z/3< \Gamma_3^2$ and the fixed point data for $\mathbb Z/3$ is $(1,2|2,2,2)$ or $(1,2|1,1,1)$.

(iii)(b)\qua For $\Gamma_3^2$, $1\to K_5\longrightarrow N(\mathbb Z/3)/\mathbb
Z/3 \longrightarrow {\sum_2}\to 1$, where $\mathbb Z/3< \Gamma_3^2$ and the fixed point data for $\mathbb Z/3$ is $(1,1|2,1,1)$.

(iii)(c)\qua For $\Gamma_3^2$, $N(\mathbb Z/3)/\mathbb
Z/3$ is a finite index subgroup of $\tilde \Gamma_1^2$, where $\mathbb Z/3< \Gamma_3^2$ and the fixed point data for $\mathbb Z/3$ is $(1,2|)$.

(iv)(a)\qua For $\Gamma_3^1$, $1\to K_5\longrightarrow N(\mathbb Z/3)/\mathbb
Z/3 \longrightarrow {\sum_4}\to 1$, where $\mathbb Z/3< \Gamma_3^1$ and the fixed point data for $\mathbb Z/3$ is $(1|2,2,2,2)$.

(iv)(b)\qua For $\Gamma_3^1$, $1\to K_5\longrightarrow N(\mathbb Z/3)/\mathbb
Z/3 \longrightarrow {\sum_3}\to 1$, where $\mathbb Z/3< \Gamma_3^1$ and the fixed point data for $\mathbb Z/3$ is $(1|2,1,1,1)$.

(iv)(c)\qua For $\Gamma_3^1$, $N(\mathbb Z/3)/\mathbb
Z/3$ is a finite index subgroup of $\tilde \Gamma_1^2$ , where $\mathbb Z/3\in \Gamma_3^1$ and the fixed point data for $\mathbb Z/3$ is $(1|2)$.

Case (2)(i):\qua For $\Gamma_3^4$, $1\to K_5\longrightarrow N(\mathbb Z/3)/\mathbb Z/3 \longrightarrow {\sum_1}\to 1$. It is easy to see by Cohen [2] that  
$$ H^i(N(\mathbb Z/3)/\mathbb Z/3,\mathbb Z)\cong H^i(K_5,\mathbb Z)\cong \left\{
\aligned \mathbb Z  \ \ \ \ \ \ \ \ i=0\\
5\mathbb Z  \ \ \ \ \ \ \ \ i=1\\
6\mathbb Z  \ \ \ \ \ \ \ \ i=2\\
0 \ \ \ \ \ \ \ \ i\ge 3
\endaligned
\right\}.
$$
We also have a short exact sequence $1\to \mathbb Z/3\longrightarrow N(\mathbb Z/3) \longrightarrow N(\mathbb Z/3)/\mathbb Z/3\to 1$.

Its associated spectral sequence collapses, giving
$$\hat H^i(N(\mathbb Z/3),\mathbb Z)_{(3)}=\left\{
\aligned 7\mathbb Z/3 \ \ or \ \ (5 \mathbb Z/3)\oplus \mathbb Z/9   \ \ \ \ \ \ \ \ i=0~mod(2)\\
5\mathbb Z/3  \ \ \ \ \ \ \ \ i=1~mod(2)
\endaligned
\right\}.
$$

Since there are five different possible types of fixed point data, 

namely $(1,2,2,2|2),$ or $(1,2,1,1|1)$, or $(1,1,2,1|1)$, 

or $(1,1,1,2|1)$, or $(1,1,1,1|2)$, 

so $\hat H^i(\Gamma_3^4,\mathbb Z)_{(3)}$

$=35\mathbb Z/3, \ \ or \ \ (33\mathbb Z/3)\oplus\mathbb Z/9,\ \ or\ \ (31\mathbb Z/3)\oplus(2\mathbb Z/9),$

$\ \ or \ \ (29\mathbb Z/3)\oplus(3\mathbb Z/9),  \ \ or \ \ (27\mathbb Z/3)\oplus(4\mathbb Z/9),$

$\ \ or \ \ (25\mathbb Z/3)\oplus(5\mathbb Z/9),\ \ \ \ \ \ \ \ i=0~mod(2)$

$\hat H^i(\Gamma_3^4,\mathbb Z)_{(3)}=25\mathbb Z/3  \ \ \ \ \ \ \ \ i=1~mod(2)$.

It is the same for the case of $\Gamma_3^5$.

Case (2)(ii)(b):\qua For $\Gamma_3^3$, if the fixed point data is $(1,1,1|1,2)$, there is a short exact sequence: $1\to K_5\longrightarrow N(\mathbb Z/3)/\mathbb Z/3 \longrightarrow {\sum_1}\to 1$.
 By the same argument as in case (2)(i), we have
$$\hat H^i(N(\mathbb Z/3),\mathbb Z)_{(3)}=\left\{
\aligned 7\mathbb Z/3 \ \ or \ \ \mathbb Z/9\oplus (5 \mathbb Z/3)    \ \ \ \ \ \ \ \ i=0~mod(2)\\
5\mathbb Z/3  \ \ \ \ \ \ \ \ i=1~mod(2)
\endaligned
\right\}.
$$

Case (2)(ii)(a):\qua For $\Gamma_3^3$, if the fixed point data is $(1,2,2|2,2),$ or $(1,2,1|1,1),$ or $(1,1,2|1,1)$, the short exact sequence is :

$1\to K_5\longrightarrow N(\mathbb Z/3)/\mathbb
Z/3 \longrightarrow {\sum_2}\to 1$, where $\mathbb Z/3< \Gamma_3^3.$

The Serre spectral sequence takes the form 

$E^{ij}_2=H^i(\Sigma_2,H^j(K_5,F_3))\Longrightarrow H^{i+j}(N(\mathbb Z/3)/\mathbb Z/3,F_3)$. 
Since 2 and 3 are relatively prime, as before, we only need to consider $i=0.$ 
$E^{0j}_2=H^j(K_5,F_3)^{\Sigma_2}$. This is very similar to the case in the previous theorem, however, we have 4 fixed points in the previous theorem, here the case is more complicated. We assume that the five 
fixed points of the $\mathbb Z/3$ action are $P_1,P_2,P_3,P_4,P_5$. Because the fixed point data are $(1,2,2|2,2),$ or $(1,2,1|1,1),$ or $(1,1,2|1,1)$, the elements of the normalizer $N(\mathbb Z/3)$ fix three points and permute the other two. Without loss of generality, $P_1,P_2,P_3$ are fixed and $P_4, P_5$ are permuted, so, $\Sigma_2$ 
is generated by $<x>\ \ =\ \ <(45)>$. By Cohen in [2], $H^1(K_5,F_3)$ is generated by five degree-one generators $\{B_{42}, B_{43}, B_{52}, B_{53},  B_{54}\}$. Also, the actions of $(45)$ on $H^1(K_5,\mathbb Z)$ are as following: $(45)B_{42}=-B_{52}$, $(45)B_{43}=B_{53}$, $(45)B_{52}=B_{42},$ $(45)B_{53}=B_{43},$ and $(45)B_{54}=B_{54}+B_{53}+B_{52}-B_{43}-B_{42}$. 
So, the invariant is generated by $<B_{42}+B_{52},\ \ B_{43}+B_{53},\ \ 2B_{42}+2B_{43}+B_{54}>$.
Thus $$H^1(K_5,F_3)^{<x>}=
<B_{42}+B_{52},\ \ B_{43}+B_{53},\ \ 2B_{42}+2B_{43}+B_{54}>\cong F_3\oplus F_3\oplus F_3.$$ Also by Cohen [2], $H^2(K_5,\mathbb Z)$ is generated by six degree-two generators 

$\{B_{42}B_{52},\ \  B_{42}B_{53},\ \ B_{42}B_{54},\ \ B_{43}B_{52},\ \  B_{43}B_{53},\ \ B_{43}B_{54}\}$. The actions of $(45)$ on $H^2(K_5,F_3)$ are induced by the actions of $(45)$ on $H^1(K_5,F_3).$ (We omit the details, it can be found in [2] or [4]) So, the invariant is generated by $<2B_{42}B_{52}+B_{42}B_{54}+B_{43}B_{52},\ \ B_{42}B_{53}+2B_{43}B_{52},\ \ 2B_{43}B_{53}+B_{43}B_{54}>$. Thus $$H^2(K_5,F_3)^{<x>}$$
$$=<2B_{42}B_{52}+B_{42}B_{54}+B_{43}B_{52},\ \ B_{42}B_{53}+2B_{43}B_{52},\ \ 2B_{43}B_{53}+B_{43}B_{54}>$$
$$\cong F_3\oplus F_3\oplus F_3.$$

The spectral sequence collapses. Therefore,  
$$ H^i(N(\mathbb Z/3)/\mathbb Z/3,F_3)=\left\{
\aligned F_3  \ \ \ \ \ \ \ \ i=0\\
3F_3  \ \ \ \ \ \ \ \ i=1\\
3F_3  \ \ \ \ \ \ \ \ i=2\\
0 \ \ \ \ \ \ \ \ i\ge 3
\endaligned
\right\}.
$$

Together with the Serre spectral sequence associated to the short exact sequence $1\to \mathbb Z/3\longrightarrow N(\mathbb Z/3) \longrightarrow N(\mathbb Z/3)/\mathbb Z/3\to 1$, this allows us to find 
$$ H^i(N(\mathbb Z/3),\mathbb Z)_{(3)}=\left\{
\aligned 4\mathbb Z/3 \ \ or \ \ (2 \mathbb Z/3)\oplus \mathbb Z/9 \ \ \ \ \ \ \ \ i=0~mod(2)\\
3\mathbb Z/3  \ \ \ \ \ \ \ \ i=1~mod(2)
\endaligned
\right\}.
$$

We put case (2)(ii)(a) and case (2)(ii)(b) together to get  

$\hat H^i(\Gamma_3^3,\mathbb Z)_{(3)}=19\mathbb Z/3 \ \ or \ \ (17\mathbb Z/3)\oplus\mathbb Z/9\ \ or\ \ (15\mathbb Z/3)\oplus(2\mathbb Z/9)\ \ or \ \ (13\mathbb Z/3)\oplus(3\mathbb Z/9)\ \ or\ \ (11\mathbb Z/3)\oplus(4\mathbb Z/9)  \ \ \ \ \ \ \ \ i=0~mod(2)$

$\hat H^i(\Gamma_3^3,\mathbb Z)_{(3)}=14\mathbb Z/3  \ \ \ \ \ \ \ \ i=1~mod(2)$.

Case (2)(iii)(b):\qua For $\Gamma_3^2$, if the fixed point data is $(1,1|2,1,1)$, there is a short exact sequence: $1\to K_5\longrightarrow N(\mathbb Z/3)/\mathbb Z/3 \longrightarrow {\sum_2}\to 1$.
 As in case (2)(ii)(a), $$\hat H^i(N(\mathbb Z/3),\mathbb Z)_{(3)}=\left\{
\aligned 4\mathbb Z/3 \ \ or \ \ \mathbb Z/9\oplus (2 \mathbb Z/3)    \ \ \ \ \ \ \ \ i=0~mod(2)\\
3\mathbb Z/3  \ \ \ \ \ \ \ \ i=1~mod(2)
\endaligned
\right\}.
$$

Case (2)(iii)(a):\qua For $\Gamma_3^2$, if the fixed point data is $(1,2|2,2,2)$ or $(1,2|1,1,1)$. There is a short exact sequence is : $1\to K_5\longrightarrow N(\mathbb Z/3)/\mathbb
Z/3 \longrightarrow {\sum_3}\to 1$, where $\mathbb Z/3\in \Gamma_3^2.$

The Serre spectral sequence takes the form 

$E^{ij}_2\cong H^i(\Sigma_3,H^j(K_5,F_3))\Longrightarrow H^{i+j}(N(\mathbb Z/3)/\mathbb Z/3,F_3)$. This case is much more complicated than before. Let us assume that the five 
fixed points of the $\mathbb Z/3$ action are $P_1,P_2,P_3,P_4,P_5$ and the elements of the normalizer $N(\mathbb Z/3)$ fix two points and permute the other three because of the fixed point data type $(1,2|2,2,2)$ or $(1,2|1,1,1)$. Without loss of generality, we assume that $P_1,P_2$ are fixed and $P_3,P_4, P_5$ are permuted, so, $\sum_3$ 
is generated by $<x,y>\ \ =\ \ <(34),(45)>$. Cohen in [2] shows that, $H^1(K_5,\mathbb Z)$ is generated by five degree-one generators $\{B_{42}, B_{43}, B_{52}, B_{53},  B_{54}\}$, and the actions are: 
$$(345)B_{42}=(34)(45)B_{42}=B_{52}-B_{42},$$ 
$$(345)B_{43}=B_{53}+B_{54},$$ 
$$(345)B_{52}=-B_{42},$$ 
$$(345)B_{53}=B_{42}+B_{43},$$ and 
$$(345)B_{54}=B_{53}+B_{52}-B_{43}-B_{42}.$$ 
So the invariant is generated by $<B_{42}+B_{52},\ \ B_{43}+B_{52}+2B_{53}+B_{54}>$. 

We can get $$H^1(K_5,F_3)^{<345>}=<B_{42}+B_{52},\ \ B_{43}+B_{52}+2B_{53}+B_{54}>\cong F_3\oplus F_3. $$

We also get $$H^1(K_5,F_3)_{<345>}=<\bar B_{42}, \bar B_{43}>\cong F_3\oplus F_3.$$ 

Consider the norm map $N:H^1(K_5,F_3)_{<345>}\rightarrow H^1(K_4,F_3)^{<345>}$. It is easy to verify that $N(\bar B_{43})=B_{43}+B_{52}+2B_{53}+B_{54}$(the generator in $H^1(K_4,F_3)^{<345>}$), and $N(\bar B_{42})=0.$ 
Therefore,
$$H^i(<345>,H^1(K_5,F_3))=coker\ \ N=F_3, \ \ if\ \ i=0~mod(2)\ \ (i>0),$$
$$H^i(<345>,H^1(K_5,F_3))=ker\ \ N=F_3, \ \ if\ \ i=1~mod(2),$$
$$H^0(<345>,H^1(K_5,F_3))=H^1(K_5,F_3)^{<x>}=F_3\oplus F_3.$$

Also, by Cohen [2], $H^2(K_5,\mathbb Z)$ is generated by six degree-two generators
 
$\{B_{42}B_{52},\ \ B_{42}B_{53},\ \ B_{42}B_{54},\ \ B_{43}B_{52},\ \  B_{43}B_{53},\ \  B_{43}B_{54}\}$. It is not hard to find that the invariant is generated by $\{B_{42}B_{52},\ \  B_{42}B_{54}+B_{43}B_{52}\}$, i.e., 
$$H^2(K_5,F_3)^{<345>}=<B_{42}+B_{52},\ \ B_{43}+B_{52}+2B_{53}+B_{54}>.$$ We also get $$H^2(K_5,F_3)_{<345>}=<\bar {B_{42}B_{53}} ,\bar {B_{43}B_{53}} >.$$ (The explicit action of $<345>$ on degree-two generators can be found in [2] or [4].)
Consider the norm map $N:H^2(K_5,F_3)_{<345>}\rightarrow H^2(K_5,F_3)^{<345>}$. One can verify that this is an isomorphism.  
Therefore, 
$$H^i(<345>,H^2(K_5,F_3))=coker\ \ N=0, \ \ if\ \ i=0~mod(2)\ \ (i>0),$$
$$H^i(<345>,H^2(K_5,F_3))=ker\ \ N=0, \ \ if\ \ i=1~mod(2),$$
$$H^0(<345>,H^2(K_5,F_3))=H^2(K_5,F_3)^{<345>}=F_3\oplus F_3.$$

Also,

$$H^i(<345>, H^j(K_5,F_3))=0, \ \ for\ \ j>2, $$
$$H^i(<345>, H^0(K_5,F_3))=H^i(<345>, F_3)=F_3, \ \ for\ \ i\ge 0$$

Note that our answers are similar to the results in [7] Lemma 4.1 and Lemma 4.2, where Xia handles unpunctured mapping class groups. Following the same calculation as in [7], we get result as Proposition 5.2 in [7].  

$$\hat H^i(N(\mathbb Z/3),\mathbb Z)_{(3)}=\left\{
\aligned \mathbb Z/3\oplus \mathbb Z/9   \ \ \ \ \ \ \ \ i=0~mod(2)\\
\mathbb Z/3  \ \ \ \ \ \ \ \ i=1~mod(2)
\endaligned
\right\}.
$$

We need to specify why we can use Xia's result from [7].

This is because:

(1)\qua Any imbedding from $\sum_3\to \sum_4$ induces an isomorphism 

$\hat H^q(\sum_4, A)_{(3)}\cong \hat H^q(\sum_3,A)_{(3)}.$ 

For $N(\mathbb Z/3)$ are same in $\sum_3$ and $\sum_4$, i.e., the ``stable" cohomology classes are the same for $\sum_3$ and $\sum_4$.

(2)\qua Even though $\Gamma_3$ has period 4, Xia proved that $1\to\pi_2\to
N(\pi_2)\to N(\pi_2)/\pi_2\to 1$ is a central extension. $N(\pi_2)$
acts on $\pi_2$ trivially, which allows him to apply the spectral
sequence argument. In our case, the period is 2 (see [4]), which
guarantees that $N(\mathbb Z/3)$ acts on $\mathbb Z/3$ trivially.
 
Case (2)(iii)(c):\qua For $\Gamma_3^2$, if the fixed point data is $(1,2|),$ $N(\mathbb Z/3)/\mathbb Z/3$ is a finite index subgroup of $\tilde \Gamma_1^2$. By Proposition 3.6-3.8 in [7], we have 
$$ H^i(N(\mathbb Z/3)/\mathbb Z/3,F_3)=0\ \ if\ \ i>0$$
$$ H^0(N(\mathbb Z/3)/\mathbb Z/3,F_3)=F_3$$
In contrast to [7], we have the central extension $1\to \mathbb Z/3\to N(\mathbb Z/3)\to N(\mathbb Z/3)/\mathbb Z/3\to 1$. Note that $N(\mathbb Z/3)$ acts on $\mathbb Z/3$ trivially. We have that 
$$\hat H^i(N(\mathbb Z/3),\mathbb Z)_{(3)}=\left\{
\aligned \mathbb Z/3   \ \ \ \ \ \ \ \ i=0~mod(2)\\
0  \ \ \ \ \ \ \ \ i=1~mod(2)
\endaligned
\right\}.
$$

Indeed, any extension $1\to \mathbb Z/p\to N(\mathbb Z/p)\to N(\mathbb Z/p)/\mathbb Z/p\to 1$ in $\Gamma_g^i$ 
is central or not can be detected by the fixed point data.

In [7], $\mathbb Z/3=\ \ <\alpha>$ and $\delta(\alpha)=(1,2)$. We have $\delta(\alpha^2)=(2,1).$ As an unordered tuple, $(1,2)=(2,1)$. Thus $\alpha$ is conjugate to $\alpha^2$, and the normalizer of $\mathbb Z/3$ is not the centralizer of $\mathbb Z/3$.

In our case, $\mathbb Z/3=\ \ <\alpha>$ and $\delta(\alpha)=(1,2|)$, and $\delta(\alpha^2)=(2,1|).$ As an ordered tuple, $(1,2|)\neq (2,1|)$. So $\alpha$ is $\bf{not}$   conjugate to $\alpha^2$, and the normalizer of $\mathbb Z/3$ is the centralizer of $\mathbb Z/3$.

We put Case (2)(iii)(a),(b),(c) together to get  

$\hat H^i(\Gamma_3^2,\mathbb Z)_{(3)}=(6\mathbb Z/3)\oplus\mathbb Z/9\ \ or\ \ (4\mathbb Z/3)\oplus(2\mathbb Z/9)  \ \ \ \ \ \ \ \ i=0~mod(2)$

$\hat H^i(\Gamma_3^2,\mathbb Z)_{(3)}=4\mathbb Z/3  \ \ \ \ \ \ \ \ i=1~mod(2)$.

Case (2)(iv)(a):\qua For $\Gamma_3^1$, if the fixed point data for $\mathbb Z/3$ is $(1|2,2,2,2),$ we have an extension $1\to K_5\longrightarrow N(\mathbb Z/3)/\mathbb
Z/3 \longrightarrow {\sum_4}\to 1.$ By Proposition 5.2 in [7], we have 
$$
\hat H^i(N(\mathbb Z/3),\mathbb Z)_{(3)}=\left\{
\aligned \mathbb Z/3\oplus \mathbb Z/9   \ \ \ \ \ \ \ \ i=0~mod(2)\\
\mathbb Z/3  \ \ \ \ \ \ \ \ i=1~mod(2)
\endaligned
\right\}.
$$

Case (2)(iv)(b):\qua For $\Gamma_3^1$, if the fixed point data for $\mathbb Z/3$ is $(1|2,1,1,1)$, we have an extension $1\to K_5\longrightarrow N(\mathbb Z/3)/\mathbb
Z/3 \longrightarrow {\sum_3}\to 1.$ Then as in Case (2)(iii)(a),  
$$
\hat H^i(N(\mathbb Z/3),\mathbb Z)_{(3)}=\left\{
\aligned \mathbb Z/3\oplus \mathbb Z/9   \ \ \ \ \ \ \ \ i=0~mod(2)\\
\mathbb Z/3  \ \ \ \ \ \ \ \ i=1~mod(2)
\endaligned
\right\}.
$$

Case (2)(iv)(c):\qua For $\Gamma_3^1$, if the fixed point data for $\mathbb Z/3$ is $(1|2),$ $N(\mathbb Z/3)/\mathbb Z/3$ is a finite index subgroup of $\tilde \Gamma_1^2$. Anologously to  Case (2)(iii)(c), 
$$
\hat H^i(N(\mathbb Z/3),\mathbb Z)_{(3)}=\left\{
\aligned \mathbb Z/3 \ \ \ \ \ \ \ \ i=0~mod(2)\\
0  \ \ \ \ \ \ \ \ i=1~mod(2)
\endaligned
\right\}.
$$

Case (2)(iv)(a),(b),(c) together, imply 

$\hat H^i(\Gamma_3^1,\mathbb Z)_{(3)}=(3\mathbb Z/3)\oplus(2\mathbb Z/9)  \ \ \ \ \ \ \ \ i=0~mod(2)$

$\hat H^i(\Gamma_3^1,\mathbb Z)_{(3)}=2\mathbb Z/3  \ \ \ \ \ \ \ \ i=1~mod(2)$.

\end{proof}

\begin{theorem}
\label{h}
For $p>3,$
\begin{enumerate}\romanIndex 
\item 
$$
\hat H^*(\Gamma_p^1,\mathbb Z)_{(p)}=\left\{
\aligned \mathbb Z/p  \ \ \ \ \ \ \ \ i=0~mod(2)\\
0  \ \ \ \ \ \ \ \ i=1~mod(2)
\endaligned
\right\}.
$$
\item
$$
\hat H^*(\Gamma_p^2,\mathbb Z)_{(p)}=\left\{
\aligned \mathbb Z/p  \ \ \ \ \ \ \ \ i=0~mod(2)\\
0  \ \ \ \ \ \ \ \ i=1~mod(2)
\endaligned
\right\}.
$$
\item
$$
\hat H^*(\Gamma_p^i,\mathbb Z)_{(p)}=0, for\ \ i\ge 3.
$$
\end{enumerate}
\end{theorem}

\begin{proof}

Plug $g=p$ into the Riemann Hurwitz equation $2p-2=p(2h-2)+t(p-1).$ For $p>3,$ the only non-negative integer solution is $(h,t)=(1,2)$. So there is no $p$-torsion in $\Gamma_p^i$ for $i\ge 3$. For $\Gamma_p^1$, by Theorem 2.7 in [4], the $\mathbb Z/p$ action on $S_p$ has $2$ fixed points and the quotient space is $S_1$.
Following the similar arguments as in [7], we have:
$$ H^i(N(\mathbb Z/p)/\mathbb Z/p,F_p)=0\ \ if\ \ i>0$$
$$ H^0(N(\mathbb Z/p)/\mathbb Z/p,F_p)=F_p$$
As before, we have the central extension $1\to \mathbb Z/p\to N(\mathbb Z/p)\to N(\mathbb Z/p)/\mathbb Z/p\to 1$. So $N(\mathbb Z/p)$ acts on $\mathbb Z/p$ trivially, and we have 
$$\hat H^i(N(\mathbb Z/p),\mathbb Z)_{(p)}=\left\{
\aligned \mathbb Z/p   \ \ \ \ \ \ \ \ i=0~mod(2)\\
0  \ \ \ \ \ \ \ \ i=1~mod(2)
\endaligned
\right\}.
$$

Since there is only one type of fixed point data $(1|p-1)$, we have one conjugacy class of $\mathbb Z/p$ in $\Gamma_p^1$. Therefore,
$$
\hat H^*(\Gamma_p^1,\mathbb Z)_{(p)}=\left\{
\aligned \mathbb Z/p  \ \ \ \ \ \ \ \ i=0~mod(2)\\
0  \ \ \ \ \ \ \ \ i=1~mod(2)
\endaligned
\right\}.
$$

The proof for $\Gamma_p^2,~p>3$, is similar.

\end{proof}

The author thanks my colleague Ethan Berkove and the referee for their suggestions.

\Addresses\recd
 
\end{document}